\documentclass[a4paper,reqno]{amsart}
\usepackage{amscd, amssymb, pb-diagram, amsmath}
\usepackage{amsthm}
\usepackage[bookmarks=true]{hyperref} 
\usepackage{slashed}

\usepackage{tikz}

\usepackage{tikz-cd}
\usepackage[margin=3.7cm]{geometry}
\usepackage{color}
\usepackage{tensor}
\usetikzlibrary{arrows,matrix,backgrounds,decorations.markings}
\numberwithin{equation}{section}

\usepackage[all]{xy}

\DeclareFontFamily{OT1}{pzc}{}
\DeclareFontShape{OT1}{pzc}{m}{it}{<-> s * [1.2] pzcmi7t}{}
\DeclareMathAlphabet{\mathpzc}{OT1}{pzc}{m}{it}

\def\End{\operatorname{End}}

\def\Dom{\operatorname{Dom}}

\def\dim{\operatorname{dim}}

\def\rd{\operatorname{d}\!}

\def\Tr{\operatorname{Tr}}

\def\sup{\operatorname{sup}}

\def\C{\mathbb{C}}

\def\R{\mathbb{R}}
\def\N{\mathbb{N}}
\def\Z{\mathbb{Z}}

\newtheorem{thm*}{Theorem}
\newtheorem{thm}{Theorem}[section]

\newtheorem{cor}[thm]{Corollary}
\newtheorem{cor*}[thm*]{Corollary}

\newtheorem{lemma}[thm]{Lemma}
\newtheorem{prop}[thm]{Proposition}

\newtheorem{problem*}[thm*]{Problem}

\theoremstyle{definition}
\newtheorem{definition}[thm]{Definition}
\newtheorem{definition*}[thm*]{Definition}

\newtheorem{problem}[thm]{Problem}

\theoremstyle{remark}
\newtheorem{remark}[thm]{Remark}
\newtheorem{remark*}[thm*]{Remark}

\newtheorem{example}[thm]{Example}

\begin{document}

\date{\today}

\title[Finite-summability in noncommutative geometry]{Finite-summability in the noncommutative geometry of crossed products}

\author{Robin J. Deeley, Magnus Goffeng, Levi Lorenzo}

\address[R. J. Deeley]{Department of Mathematics, University of Colorado Boulder Campus Box 395, Boulder, CO 80309-0395, USA}
\email{robin.deeley@colorado.edu}

\address[M. Goffeng]{Centre for Mathematical Sciences, Lund University, Sweden}
\email{magnus.goffeng@math.lth.se}

\address[L. Lorenzo]{Department of Mathematics, Fort Lewis College, 1000 Rim Drive, Durango, Colorado 81301-3999, USA}
\email{llorenzo@fortlewis.edu}

\thanks{RJD was partially supported by NSF Grants DMS 2000057, DMS 2247424 and Simons Foundation Gift MP-TSM-00002896. MG was supported by the Swedish Research Council Grant VR 2025-03923. LL was also partially supported by NSF Grants DMS 2000057 and DMS 2247424 while being the PhD student of RJD}

\maketitle

\begin{abstract}
We study finite-summability of unbounded representatives of $K$-homology classes for crossed-product $C^*$-algebras associated with actions on compact spaces. Using Kasparov's Dirac--dual-Dirac method, we isolate geometric conditions on the $\gamma$-element under which elements of $K$-homology can be transferred to a crossed product associated with a proper action and represented by finitely summable unbounded cycles defined from Dirac-type operators with an unbounded potential. We relate finite-summability of these operators to integrability properties of the potential by spectral estimates for Schrödinger operators on manifolds of bounded geometry. In particular, for finitely generated nilpotent groups this leads to finitely summable, unbounded representatives of all $K$-homology classes. We also give examples showing that representability of $K$-homology classes by finitely summable unbounded cycles is not invariant under $KK$-equivalence.
\end{abstract}

\section{Introduction}

In this paper we study finite-summability in noncommutative geometry. In Connes' spectral noncommutative geometry \cite{redbook}, geometry is encoded through spectral triples allowing for techniques of differential geometry and tools from global analysis to be extended to general ``noncommutative spaces'' -- geometric objects corresponding to $C^*$-algebras hinted at by the Gelfand-Naimark theorem (but not existing in classical sense). The relevant notion of dimension in Connes' spectral noncommutative geometry is that of finite-summability degree, connecting to classical dimensions via the Weyl law for elliptic operators \cite{horIV}. Finite-summability plays an important role internally in noncommutative geometry in several applications as well as for defining Connes-Chern character and for the celebrated local index formula \cite{connesmoscovici}. We approach this sphere of thought by means of another important tool from noncommutative geometry, the techniques derived from the work on the Baum-Connes conjecture by Kasparov \cite{kaspinv} relying on $\gamma$-elements and special actions.

The problem motivating this work is that of determining when finite-summability is attainable. This problem is already well studied in the literature \cite{EmDuw,EN, GM15,GRU, Pusch} for special classes of $C^*$-algebras $A$. Our contribution is to show how Kasparov's method for proving the Baum-Connes conjecture via $\gamma$-elements can be used to prove finite-summability under geometric constraints on the $\gamma$-element. We here mean attainable in the rather weak sense that we ask which $K$-homology classes on a $C^*$-algebra can be represented by finitely summable cycles. In other words, we study when a spectral noncommutative geometry up to homotopy is finitely summable. The problem comes in two flavors as we can take cycles to be ordinary bounded $KK$-cycles or Baaj-Julg's unbounded cycles -- closely related to Connes' notion of a spectral triple. To study when $K$-homology classes can be represented by finitely summable bounded cycles is a delicate problem with only one known counterexample \cite{GM15}, even if a seminal result of Puschnigg \cite{Pusch} indicates that no non-trivial $K$-homology class on the reduced group $C^*$-algebra of a higher rank lattice can be represented by a finitely summable bounded cycle.

We shall study the problem of when all $K$-homology classes of a $C^*$-algebra $A$ are representably by finitely summable unbounded cycles. In more colloquial terms, we ask when a noncommutative geometry can be described by a spectral triple encoding a finite-dimensional system. Since any finitely summable unbounded cycle produces a finitely summable bounded cycle in the same $K$-homology class \cite{baajjulg,sww}, it is a stronger property to ask for all $K$-homology classes in a $C^*$-algebra to be representable by finitely summable unbounded cycles. There are severe obstructions on a $C^*$-algebra to admit finitely summable unbounded cycles. In fact Connes \cite{ConnesETDS} proved that any finitely summable unbounded cycle on a $C^*$-algebra gives rise to a tracial state on it, a result later relaxed in \cite{GRU} to also include super-logarithmic growth of eigenvalues. Therefore $C^*$-algebras without tracial states, e.g. purely infinite $C^*$-algebras, admit no non-trivial finitely summable unbounded cycles and only the trivial $K$-homology class is representable by a finitely summable unbounded cycle.

Our contribution to the understanding of the problem discussed in the preceding paragraph is to show that it is indeed hard in the sense that it is not invariant under $KK$-isomorphism, and we also show that it has a positive solution for $A=C(X)\rtimes \Gamma$ whenever $X$ is a finite CW-complex with an action of a finitely generated nilpotent group $\Gamma$. The former negative result is found from a particularly nice $KK$-isomorphism coming from the Dirac element of a surface group. The latter result, despite being somewhat expected due to its finite-dimensional nature, has not been found by the authors in the literature. Our technique breaks down at the spectral level when $\Gamma$ has super-polynomial growth, for instance for hyperbolic groups or higher rank lattices. We also demonstrate how our results apply to examples such as flat manifolds and the noncomutative torus. Throughout the paper, crossed products will refer to the full crossed product. These results are discussed in more detail in the next section.

\subsection*{Acknowledgements} The authors wish to thank Rupert Frank, Ada Masters, Ryszard Nest and Grigori Rozenblum for discussions. 

\subsection*{Usage of AI} 

In the writing of this paper, ChatGPT 5.6 was used for proofreading and editorial purposes, as well as helping with the proof of Lemma \ref{schattenbounds}.

\section{Main results}
\label{subsecintro}

Our main results concern finite-summability of unbounded representatives of $K$-homology classes on algebras of the form $C(X)\rtimes \Gamma$, where $\Gamma$ is a discrete group acting on a finite CW-complex $X$. We often impose conditions on the $\gamma$-element of $\Gamma$. Recall from \cite{tu00} that the $\gamma$-element of a discrete group $\Gamma$, when it exists, is a distinguished element $\gamma\in KK^\Gamma(\C,\C)$ satisfying the following conditions:
\begin{enumerate}
\item[A.] For any finite subgroup $F\subseteq \Gamma$, $\mathrm{Res}^\Gamma_F(\gamma)=1_F\in R(F)\cong KK^F(\C,\C)$.
\item[B.] There is a proper $\Gamma-C^*$-algebra $\mathcal{P}$, and elements $\alpha\in KK^\Gamma(\mathcal{P},\C)$ and $\beta\in KK^\Gamma(\C,\mathcal{P})$ such that 
$$\gamma=\beta\otimes_\mathcal{P} \alpha\in KK^\Gamma(\C,\C).$$
\end{enumerate}
For an alternate approach, see also \cite{nishisplit}. For large classes of groups there exists a $\gamma$-element \cite{kaspinv}, for instance for closed subgroups of almost connected groups. Recall that the $\gamma$-element is unique when it exists and that the mere existence of a $\gamma$-element for $\Gamma$ implies split injectivity of both the full and reduced Baum-Connes assembly map for $\Gamma$ with arbitrary coefficients \cite{tu00}. We remark here that $\gamma=1$ is equivalent to the strong Baum-Connes conjecture, cf. \cite{meyernest}, and is known for a-T-menable groups by work of Higson-Kasparov \cite{higsonkasp}. The Baum-Connes conjecture with coefficients is equivalent to $\gamma$ acting as $1$ on $K$-theory of crossed products by $\Gamma$, but in general the Baum-Connes conjecture with coefficients does not say anything about how $\gamma$ acts on $K$-homology of crossed products by $\Gamma$. We consider the following further conditions on a $\gamma$-element.

\begin{definition*}
\label{geomegamadma}
Let $\Gamma$ be a countable group with $\gamma$-element $\gamma\in KK^\Gamma(\C,\C)$. 
\begin{itemize}
\item We say that the $\gamma$-element has a geometric model if in the factorization $\gamma=\beta\otimes_\mathcal{P} \alpha$ we can take $\mathcal{P}=C_0(Z,\mathcal{A})$, where $Z$ is a proper, $\Gamma$-finite $\Gamma$-CW-complex $Z$ and $\mathcal{A}\to Z$ is a finite rank Azumaya bundle, and $\alpha\in KK^\Gamma(\mathcal{P},\C)$ and $\beta\in KK^\Gamma(\C,\mathcal{P})$ to satisfy that 
$$1_\mathcal{P}=\alpha\otimes_\C \beta\in KK^\Gamma(\mathcal{P},\mathcal{P}).$$
\item Let $m>0$. We say that $\gamma$ admits a geometric $\beta$-factor of order $m$ if $\gamma$ additionally to being geometric satisfies that $Z$ can be chosen to be a proper, connected, cocompact $\Gamma$-manifold with an isometric action on a Riemannian metric, and $\beta$ can be represented by an order $m$ unbounded $\Gamma$-equivariant cycle of the form $(C_0(Z,\mathcal{E}),V)$ where 
\begin{itemize}
\item $\mathcal{E}\to Z$ is a hermitean vector bundle of right $\mathcal{A}$-modules; and 
\item $V$ is a locally Lipschitz endomorphism valued section such that $V^2$ is scalar valued and for some $\alpha>0$, $V$ satisfies the bound $|\nabla V|\lesssim (1+|V|)^{1-\frac{1}{\alpha}}$. 
\end{itemize}
\end{itemize}
The notion of order $m$ unbounded $\Gamma$-equivariant cycles is recalled below in Definition \ref{adljnadjlknad}, this notion encodes that $V$ need not be invariant up to bounded terms but only invariant up to terms relatively bounded by a small power of $V$. We tacitly use the notations of this definition to indicate the geometric data in the factors associated with a geometric $\gamma$-element or a $\gamma$-element with a geometric $\beta$-factor, even if the geometric data is not uniquely determined.
\end{definition*}

The order $m$ measures the equivariance defect of the unbounded cycle, whereas the exponent $\alpha$ controls the growth of the derivative of the potential.
They are logically independent. For finitely generated nilpotent groups, see Lemma \ref{lemaodknaldn},  they can be chosen equal. We note that if $\gamma$ has a geometric $\beta$-factor of order $m$ represented by $(C_0(Z,\mathcal{E}),V)$ then for $s\in (0,1)$ and $\hat{V}:=V(1+V^2)^{-(m-1)/2m}$ the collection $(C_0(Z,\mathcal{E}),\hat{V})$ is a $\Gamma$-equivariant cycle of order $1$ representing $\beta$. In several examples, the geometry naturally gives rise to higher order representatives of the dual Dirac element $\beta$.

\begin{remark*}
Whenever $\Gamma$ admits a proper, cocompact $\Gamma$-manifold with non-positive sectional curvatures, \cite{kaspinv} ensures that $\Gamma$ has a $\gamma$-element with a geometric $\beta$-factor of order $1$ with $Z$ being the special $\Gamma$-manifold and $\mathcal{A}\to Z$ being the bundle of Clifford algebras. In particular, if $\Gamma$ is a discrete cocompact subgroup of a semisimple Lie group it has a $\gamma$-element with a geometric $\beta$-factor of order $1$. 

However, when $\Gamma$ is a lattice in a non-abelian nilpotent Lie group the natural models for the $\gamma$-element has a geometric $\beta$-factor of order $m> 1$. The reason is that the anisotropic geometry of a nilpotent Lie group utilized in standard constructions of the dual Dirac element $\beta$ naturally produces higher order polynomials.
\end{remark*}

\begin{remark*}
In several cases where $\gamma=1$ and the Baum-Connes conjecture is known to hold true, e.g. for general a-T-menable groups, the existing models for the $\gamma$-element are not geometric, for further references see the discussion in \cite[Introduction]{laffinv}.
\end{remark*}

\begin{thm*}
\label{mainthm1}
Let $\Gamma$ be a countable group with a geometric model for the $\gamma$-element. Then for any finite CW-complex $X$ with a $\Gamma$-action such that $X\times Z$ is a $\Gamma$-finite $\Gamma$-CW-complex, any element of $K^*(C_0(X\times Z,\mathcal{A})\rtimes \Gamma)$ can be represented by a finitely summable unbounded cycle. In particular, any element of $\gamma K^*(C(X)\rtimes \Gamma)$ is up to the $KK$-isomorphism 
\begin{align*}
\alpha:\gamma KK_*^\Gamma(C(X),\C)&=\gamma K^*(C(X)\rtimes \Gamma)\\
&\xrightarrow{\sim} KK_*^\Gamma(C_0(X\times Z, \mathcal{A}),\C)=KK_*(C_0(X\times Z, \mathcal{A})\rtimes \Gamma,\C),
\end{align*}
representable as a finitely summable unbounded cycle. Moreover, suppose that $\gamma$ has a geometric $\beta$-factor of order $m\geq 1$ and that $X\times Z$ is equipped with a $\Gamma$-equivariant spin$^c$-structure compatible with the geometric model. Then, if $\Gamma$ acts Lipschitz continuously on $X$, any element $\delta\in \gamma K^*(C(X)\rtimes \Gamma)$ can be represented as an $m$th order spectral triple of the form
$$(\mathrm{Lip}(X)\rtimes_{\rm alg}\Gamma,L^2(X\times Z;\mathcal{E}\otimes E\otimes S_{X\times Z} ),1\otimes_{\nabla_\mathcal{E}}\slashed{D}+V\otimes 1_{E\otimes S_{X\times Z}}),$$
for some hermitean vector bundle $E\to X\times Z$ with a twisted spin$^c$-Dirac operator $\slashed{D}$ on $E\otimes S_{X\times Z}$.
\end{thm*}

The reader can find the definition of an $m$:th order spectral triple in Definition \ref{defspectrip}. We note that we need to go to higher order spectral triples despite the relevant operator being first-order because the additional potential $V$ arising from the dual Dirac element could have poor equivariance properties quantified by the order $m$ of $\beta$. This forces the order up to $m\geq 1$.

\begin{remark*}
The utility of Theorem \ref{mainthm1} in index computations is similar to the philosophy explored in the third author's PhD-thesis \cite{levithesis}, see also \cite{levipaper}. Let us summarize the idea. Given an element $\delta\in \gamma K^*(C(X)\rtimes \Gamma)$ represented as 
\begin{align}
\label{kasprodforsda}
\delta=&(1\otimes \beta)\otimes_{C_0(X\times Z,\mathcal{A})} [L^2(X\times Z;E\otimes S_{X\times Z}),\slashed{D}]=\\
\nonumber
=&[L^2(X\times Z;\mathcal{E}\otimes E\otimes S_{X\times Z} ),1\otimes_{\nabla_\mathcal{E}}\slashed{D}+V\otimes 1_{E\otimes S_{X\times Z}}],
\end{align}
we can compute the index pairing of $K_*(C_0(X\times Z,\mathcal{A})\rtimes \Gamma)$ with $\delta':=[L^2(X\times Z;E\otimes S_{X\times Z}),\slashed{D}]$ using standard methods. In particular, for $y\in K_*(C(X)\rtimes \Gamma)$, we can using the associativity of the Kasparov product to write the index pairing of $y$ with $\delta$ as 
$$\langle y,\delta\rangle=\langle y,(1\otimes \beta)\otimes_{C_0(X\times Z,\mathcal{A})}\delta'\rangle=\langle y\otimes_{C(X)\rtimes \Gamma}\alpha ,\delta'\rangle,$$
and the latter index pairing is computable using finite-summability of the unbounded representative of $\delta'$.

\end{remark*}

Theorem \ref{mainthm1} will be proved below in Subsection \ref{proofof2}. We note that Theorem \ref{mainthm1} has consequences for permanence properties of finite-summability under $KK$-isomorphism. The Lie group $G=SU(n,1)$, for $n\geq 1$, has $\gamma_G=1$ by \cite{julgkasparov} so for any discrete subgroup $\Gamma$ of a finite quotient of $SU(n,1)$ we have that $\gamma=1$. If we take $n=1$ and $\Gamma=PSL(2,\Z)$ as a subgroup of $G=PSU(1,1)=PSL(2,\R)$, with its action by Möbius transforms on $S^1$, it was proven in \cite{ananth} that $C(S^1)\rtimes \Gamma$ is a Kirchberg algebra, in particular unital and purely infinite. We conclude the following corollary to Theorem \ref{mainthm1}.

\begin{cor*}
There are discrete subgroups $\Gamma$ of semisimple connected Lie groups $G$ such that for some compact $\Gamma$-manifold $X$ the $C^*$-algebra $C(X)\rtimes \Gamma$ is unital, purely infinite, and admits no non-trivial, finitely summable unbounded cycles, yet all $K$-homology classes on the $KK$-isomorphic $C^*$-algebra $C_0(X\times G/K,\mathcal{A})\rtimes \Gamma$ can be represented by finitely summable unbounded cycles.
\end{cor*}

We can also use Theorem  \ref{mainthm1} to obtain positive results for finitely generated nilpotent Lie groups. For such groups the dual-Dirac element admits the
higher-order geometric representative constructed in Lemma \ref{lemaodknaldn}, and the underlying nilpotent Lie group has polynomial volume growth. Combining this with the spectral estimates of Section \ref{alkdnaldn} yields finite-summability of the Kasparov product \eqref{kasprodforsda}.

\begin{thm*}
\label{mainthm3}
If $\Gamma$ is a finitely generated nilpotent group acting on a finite CW-complex $X$, then any element of $K^*(C(X)\rtimes \Gamma)$ is representable as a unital finitely summable unbounded cycle.
\end{thm*}

\begin{remark*}
By Corollary \ref{lakdjnadln}, polynomial growth is a necessary condition for our techniques to work.  We expect that Theorem \ref{mainthm3} extends to arbitrary finitely generated discrete group of polynomial growth showing that it is also a sufficient condition. However, our method of proof, in particular Lemma \ref{lemaodknaldn}, relies on $\Gamma$ being a finitely generated nilpotent group. 
\end{remark*}

\begin{remark*}
Theorem \ref{mainthm1} gives us information about the $\gamma$-part $\gamma K_*^\Gamma(X)$ of the $K$-homology of $C(X)\rtimes \Gamma$. We can ask what can be said along the lines of Theorem \ref{mainthm3} beyond groups of polynomial growth? Our method of proof, using cycles of the form \eqref{kasprodforsda}, will by our results in Section \ref{alkdnaldn} produce unbounded cycles that are not finitely summable whenever $Z$ has super-polynomial growth. 

Examples already show a substantial difference between bounded and unbounded finite-summability. Finitely summable bounded Fredholm modules have been constructed for Cuntz--Krieger algebras, Smale spaces, and reduced group $C^*$-algebras of hyperbolic groups \cite{EN, Ger,GM15}. By contrast, Connes' tracial obstruction excludes non-trivial finitely summable unbounded cycles on purely infinite $C^*$-algebras \cite{ConnesETDS,GRU}, while Puschnigg proves a strong obstruction to finitely summable Fredholm modules for higher-rank lattices \cite{Pusch}. These results suggest that bounded finite-summability may persist in situations in which global Dirac-type representatives are unavailable. It is therefore natural to ask whether methods from spectral geometry can be adapted to construct finitely summable bounded cycles beyond the range in which finitely summable unbounded cycles exist.
\end{remark*}

The paper is organized as follows. In Section \ref{subsecspecncg} we review spectral triples, finite-summability, and the bounded and unbounded pictures of $K$-homology. Section \ref{seckksda} recalls equivariant $KK$-theory and criteria for computing Kasparov products using unbounded cycles. In Section \ref{subsecgamma} we review $\gamma$-elements and the geometric models arising from special manifolds. Section \ref{alkdnaldn} studies Dirac operators with unbounded potentials on manifolds of bounded geometry and relates Schatten-class resolvent estimates to the growth of the potential. Section \ref{bdsecada} recalls geometric models for equivariant $K$-homology with Azumaya coefficients and derives finitely summable representatives for proper actions. The main results are proved in Section \ref{pfsdsdsdkskn}. In Section \ref{excxca} we study examples,  focusing on groups acting on manifolds of nonpositive sectional curvature.

\section{Preliminaries in spectral noncommutative geometry}
\label{subsecspecncg}

Before entering into the proof of our main results, let us first review and collate the work needed from the literature. As our results concern spectral noncommutative geometry, we first recall the relevant notions therefrom. There is to the knowledge of the authors no complete overview of spectral noncommutative geometry in the literature, but the reader may consult \cite{Blackadar,greenbook,yellowbook} for more elementary introductions or \cite{ConnesETDS,redbook,connesmoscovici,connesreconstr,kaspinv, laffinv,thebeastofmesland} for more advanced topics.

The fundamental object in spectral noncommutative geometry is a spectral triple. We use the notation $\mathbb{B}(\mathcal{H})$ to denote the bounded operators on a separable Hilbert space $\mathcal{H}$, normed by the operator norm. We let $\mathbb{K}(\mathcal{H})\subseteq \mathbb{B}(\mathcal{H})$ denote the norm closed ideal of compact operators on $\mathcal{H}$. Additionally, for $p\geq 1$ we write $\mathcal{L}^p(\mathcal{H})$ for the symmetrically normed ideal of compact operators 
$$\mathcal{L}^p(\mathcal{H}):=\{T\in \mathbb{B}(\mathcal{H}): \mathrm{Tr}(|T|^p)<\infty\}.$$
For more on symmetrically normed ideals, see \cite{lordzaninsuko,simonbook}.

\begin{definition}
\label{defspectrip}
A spectral triple $(\mathcal{A},\mathcal{H},D)$ of order $m>0$ consists of a separable Hilbert space $\mathcal{H}$, a $*$-subalgebra $\mathcal{A}\subseteq \mathbb{B}(\mathcal{H})$, and an unbounded, self-adjoint operator $D$ on $\mathcal{H}$ such that for any $a\in \mathcal{A}$
\begin{enumerate}
\item $a(i+D)^{-1}\in \mathbb{K}(\mathcal{H})$;
\item $a\Dom(D)\subseteq \Dom(D)$ and $[D,a](1+D^2)^{-1/2+1/2m}$ is norm bounded.
\end{enumerate}
If additionally, $p\geq 1$ satisfies that $a(1+D^2)^{-p/2}\in \mathcal{L}^1(\mathcal{H})$ for all $a\in \mathcal{A}$ we say that $(\mathcal{A},\mathcal{H},D)$ is $p$-summable. If there is a $p\geq 1$ for which $(\mathcal{A},\mathcal{H},D)$ is $p$-summable, we say that $(\mathcal{A},\mathcal{H},D)$ is finitely summable.
\end{definition}

The definition of finite-summability is well established in the unital case, and we take the definition of finite-summability in the non-unital case from \cite[Chapter 3.4]{CGRS11}. The notion of a spectral triple is flexible, and has found applications in a plethora of areas \cite{CGRS11,ConnesETDS,redbook,connesmoscovici, EmDuw,GM15,GMRtwist,GRU,kaspinv,laffinv,leschkaadsf,levipaper,thebeastofmesland}. Higher order spectral triples give rise to ordinary spectral triples (of order $1$) via logarithmic dampening \cite{GMRtwist}.

We here recall the examples relevant for the results of this paper. The reader should beware that for our applications, the spectral triples are higher order even if they arise from first-order differential operators because an additional potential arising from the dual Dirac element could have poor equivariance properties which force up the order above $1$. The prototypical example of a spectral triple of order $m$ arises from an $m$:th order, elliptic differential operator on a manifold. More precisely, we have the following proposition.

\begin{prop}
\label{adlknadlkna}
Let $M$ be a manifold with a volume density $\mu$, and assume that $E\to M$ is a Hermitean vector bundle. If $\slashed{D}:C^\infty(M,E)\to C^\infty(M,E)$ is an $m$:th order, elliptic differential operator which is essentially self-adjoint on the domain $C^\infty_c(M,E)\subseteq L^2(M,E)$ then its closure $D$ on $L^2(M,E)$ fits into a spectral triple of order $m$
\begin{equation}
\label{prororad}
(C^\infty_c(M),L^2(M,E),D).
\end{equation}
Moreover, assuming that $\Gamma$ is a discrete group of diffeomorphisms on $M$ preserving the volume density and lifting to an isometric action on $E$, then whenever the principal symbol of $\slashed{D}$ is $\Gamma$–invariant, the integrated action of $C^\infty_c(M)\rtimes_{\rm alg}\Gamma$ on $L^2(M,E)$ defines  a spectral triple 
\begin{equation}
\label{prororadcrossed}
(C^\infty_c(M)\rtimes_{\rm alg}\Gamma,L^2(M,E),D).
\end{equation}
In fact, the spectral triples \eqref{prororad} and \eqref{prororadcrossed} are $p$-summable for any $p>\dim(M)/m$.
\end{prop}

The construction of Proposition \ref{adlknadlkna} has vast extensions also to foliations, more general operators in pseudodifferential calculi and an abundance of geometric examples.

\begin{proof}
We start by noting that for $a\in C^\infty_c(M)$, we have that $a\Dom(D)\subseteq \Dom(D)$ with the operator $[D,a]$ being relatively bounded by the $(m-1)$st order elliptic operator $(1+D^2)^{1/2-1/2m}$ by elliptic regularity because $a$ preserves the core $C^\infty_c(M,E)$ of $D$ and $[D,a]$ is a compactly supported operator of order $m-1$. The fact that $a(i+D)^{-1}\in \mathbb{K}(\mathcal{H})$ follows from the Rellich theorem. In fact, by work of Seeley \cite{Seeley1967}, $a(1+D^2)^{-p/2}$ is a compactly based pseudodifferential operator of order $-pm$. Together with the Weyl law, this proves that $(C^\infty_c(M),L^2(M,E),D)$ is a $p$-summable spectral triple for any $p>\dim(M)/m$.

To prove that $(C^\infty_c(M)\rtimes_{\rm alg}\Gamma,L^2(M,E),D)$ is a $p$-summable spectral triple for any $p>\dim(M)$, it suffices to show that when $a\in C^\infty_c(M)$ and $g\in \Gamma$, we have that $ag\Dom(D)\subseteq \Dom(D)$ with the operator $[D,ag]$ being relatively bounded by $(1+D^2)^{1/2-1/2m}$. Again, this follows since $a$ preserves the core $C^\infty_c(M,E)$ of $D$ and $[\slashed{D},ag]=([\slashed{D},a]+a(D-gDg^{-1}))g$ is relatively bounded by $(1+D^2)^{1/2-1/2m}$ since $[\slashed{D},a]+a(D-gDg^{-1})$ is a compactly supported operator of order $m-1$ and $g$ acts unitarily.
\end{proof}

\begin{remark}
We remark that the assumptions made on the group action in Proposition \ref{adlknadlkna} are to ensure that the $\Gamma$-action makes $(C^\infty_c(M),L^2(M,E),D)$ into a uniformly equivariant unbounded cycle in the sense of \cite[Definition I.2.7]{adathesis}, see also \cite[Definition 8.7]{danthesis}.
\end{remark}

\begin{remark}
\label{chernoffrem}
All of the assumptions of Proposition \ref{adlknadlkna} are local conditions except that $\slashed{D}$ is essentially self-adjoint on $C^\infty_c(M,E)$. Chernoff bounds \cite{chernoffjfa} produce a sufficient condition to ensure that a first-order, elliptic, formally self-adjoint differential operator $\slashed{D}:C^\infty(M,E)\to C^\infty(M,E)$ is essentially self-adjoint on the domain $C^\infty_c(M,E)\subseteq L^2(M,E)$. The latter holds true if $M$ is a complete Riemannian manifold, with $\mu$ the Riemannian volume density, and the principal symbol $\sigma^1(\slashed{D})$ satisfies that there is a constant $C>0$ such that for $\|\sigma^1(\slashed{D})(x,\xi)\|_{\End(E_x)}\leq C|\xi|_x$ for any $(x,\xi)\in T^*M$.
\end{remark}

\begin{remark}
\label{extendtolip}
We can extend the proof of Proposition \ref{adlknadlkna} also to the Lipschitz case when $D$ has order $m=1$. Indeed, if we write $\mathrm{Lip}_c(M)$ the space of compactly supported Lipschitz functions on $M$ then under the assumptions of Proposition \ref{adlknadlkna} $(\mathrm{Lip}_c(M),L^2(M,E),D)$ forms a $p$-summable spectral triple of order $1$ for any $p>\dim(M)$. Similarly, if $\Gamma$ acts on $M$ in a way satisfying the assumptions of Proposition \ref{adlknadlkna} $(\mathrm{Lip}_c(M)\rtimes_{\rm alg}\Gamma,L^2(M,E),D)$ forms a $p$-summable spectral triple for any $p>\dim(M)$, in fact we need only require that $\Gamma$ acts by locally bi-Lipschitz homeomorphisms. The proof of these two facts extends directly in light of the factorization $\mathrm{Lip}_c(M)=C^\infty_c(M)\mathrm{Lip}_c(M)$ and elliptic regularity which ensures that $a \Dom(D)\subseteq H^1_c(M,E)$ for any $a\in C^\infty_c(M)$.
\end{remark}

\begin{example}
\label{akandaojnd}
An important special case of the construction in Proposition \ref{adlknadlkna} arises on complete Riemannian spin$^c$-manifolds. Let $M$ be a complete Riemannian spin$^c$-manifold with spinor bundle $S_M\to M$. For any hermitean vector bundle $E\to M$, we can choose a hermitean connection on $E\to M$ and a spin$^c$-connection on $S_M$ that together for a spin$^c$-Dirac operator $\slashed{D}_E:=C^\infty(M,E\otimes S_m)\to C^\infty(M,E\otimes S_M)$. Since $\slashed{D}_E$ is Dirac type, $\|\sigma^1(\slashed{D})(x,\xi)\|_{\End(E_x\otimes S_{M,x})}=|\xi|_x$ and completeness of $M$ together with Chernoff bounds (cf. Remark \ref{chernoffrem}) ensure that $\slashed{D}_E$ is essentially self-adjoint. We arrive at a $p$-summable spectral triple $(C^\infty_c(M),L^2(M,E\otimes S_M),D_E)$ of order $1$, for any $p>\dim(M)$.

If additionally $\Gamma$ acts isometrically and spin$^c$-preservingly on $M$, as well as lifting to $E$, we arrive at a $p$-summable spectral triple $(C^\infty_c(M)\rtimes_{\rm alg}\Gamma,L^2(M,E\otimes S_M),D_E)$, for any $p>\dim(M)$. We note that it is rare for infinite group actions to be isometric on a compact manifold. In fact, isometric group actions arise in many applications primarily from proper actions and in such applications, the metric is often chosen a posteriori from the characteristic of being $\Gamma$-invariant rather than considering groups acting isometrically in a fixed metric.
\end{example}

\begin{remark}
In light of the examples above, one often refers to the $*$-algebra $\mathcal{A}$ in a spectral triple as a smooth structure. A highlight of spectral noncommutative geometry is Connes' reconstruction theorem \cite{connesreconstr,renniereconst} proving that for spectral triples $(\mathcal{A},\mathcal{H},D)$ of order $1$ with $\mathcal{A}$ commutative and unital, plus an additional set of axioms, then $(\mathcal{A},\mathcal{H},D)$ is of the same form as in Example \ref{akandaojnd} with $M$ compact.
\end{remark}

For the sake of completeness, let us mention two examples of spectral triples that are truly noncommutative.

\begin{example}
For a discrete group $\Gamma$, we can define an operator-valued weight as a function $c:\Gamma\to \mathbb{B}(V)$ for a finite-dimensional Hilbert space $V$. We shall assume that $c$ is symmetric (i.e. $c(g)=c(g^{-1})$), self-adjoint (i.e. $c(\gamma)=c(\gamma)^*$), proper (i.e. $(i+c)^{-1}\in c_0(\Gamma,\mathbb{B}(V))$) and translation bounded (i.e. $\sup_g\|c(gh)-c(g)\|<\infty$ for any $h\in \Gamma$). Then we have a spectral triple 
$$(\C[\Gamma],\ell^2(\Gamma,V),D),$$
where $\C[\Gamma]$ denotes the group algebra of $\Gamma$ and $Df(g)=c(g)f(g)$ for 
$$f\in \Dom(D):=\left\{f\in \ell^2(\Gamma,V): \, \sum_{g\in \Gamma} \|c(g)f(g)\|_V^2<\infty\right\}.$$
A prototypical example often studied in the literature due to its connection with geometric group theory is $V=\C$ and $c$ being a length function \cite{ConnesETDS, GRU,adathesis}.
\end{example}

\begin{example}
Another class of examples of spectral triples arise on Cuntz-Krieger algebras \cite{GM15}, or more generally Cuntz-Pimsner algebras \cite{GMR}. Their precise form is not relevant for the purposes of this paper, however we note that they are geometric in that they carry the dynamical information and the geometry of the underlying shift spaces of finite type. What is interesting to note is that Cuntz-Krieger algebras are purely infinite when the underlying shift space is mixing, so there are not finitely summable spectral triples but rather the eigenvalues display a logarithmic growth.
\end{example}

The role of spectral triples on the ``noncommutative topology'' comes from a bounded analogue of the axioms of a spectral triple. While a spectral triple is an abstraction of first-order, elliptic operators (Proposition \ref{adlknadlkna}), bounded cycles abstracts zeroth order, elliptic operators -- this was their motivating feature as cycles for $K$-homology in Atiyah's work \cite{atiyah70}.

\begin{definition}
\label{bddcycle}
A bounded cycle on a $C^*$-algebra $A$ is a pair $(\mathcal{H},F)$ where $\mathcal{H}$ is a separable Hilbert space carrying a representation $A\to \mathbb{B}(\mathcal{H})$ and $F\in \mathbb{B}(\mathcal{H})$ is an operator satisfying for any $a\in A$ that
$$[F,a], \; a(F^*-F), \;a(F^2-1)\in \mathbb{K}(\mathcal{H}).$$
For a $p\geq 1$ and a dense $*$-subalgebra $\mathcal{A}\subseteq A$ such that for any $a\in \mathcal{A}$ we have that
$$[F,a]\in \mathcal{L}^p(\mathcal{H}),\quad \mbox{and} \quad a(F^*-F), \;a(F^2-1)\in \mathcal{L}^p(\mathcal{H}),$$
we say that $(\mathcal{H},F)$ is $p$-summable, and more precisely that $(\mathcal{H},F)$ is $p$-summable on $\mathcal{A}$. If there is a $p\geq 1$ such that $(\mathcal{H},F)$ is $p$-summable on a dense subalgebra, we say that $(\mathcal{H},F)$ is finitely summable.

When we wish to emphasize the dependence on the representation $\pi:A\to \mathbb{B}(\mathcal{H})$, we sometimes write the bounded cycle as $(\pi,\mathcal{H},F)$
\end{definition}

The connection between spectral triples and unbounded cycles come from the following result due to Baaj-Julg \cite{baajjulg}, later extended to several other cases throughout the literature. We say that $\chi\in C_b(\R)$ is a choice function if $\chi$ is an odd function and $\chi^2-1\in C_0(\R)$.

\begin{thm}[Baaj-Julg \cite{baajjulg}]
\label{baajjulgthm}
Suppose that $(\mathcal{A},\mathcal{H},D)$ is a spectral triple of order $m>0$. For any choice function $\chi$ and any representation $\pi:A\to \overline{\mathcal{A}}^{\|\cdot\|_{\mathbb{B}(\mathcal{H})}}$, the collection $(\pi,\mathcal{H},\chi(D))$ is a bounded cycle on $A$.

Moreover, if for some $p\geq1$, $(\mathcal{A},\mathcal{H},D)$ is $p$-summable and $\mathcal{A}_0\subseteq A$ is a dense $*$-subalgebra with $\pi(\mathcal{A}_0)\subseteq \mathcal{A}^2$ then $(\pi,\mathcal{H},\chi(D))$ is $p$-summable on $\mathcal{A}_0$ as soon as $\chi(x)^2-1=O(|x|^{-2})$ when $|x|\to \infty$.
\end{thm}

The proof of Theorem \ref{baajjulgthm} is found for $\chi(t)=t(1+t^2)^{-1/2}$ in \cite{baajjulg} in order $m=1$ and works the same for any $m>0$, see \cite[Theorem A.6]{GM15}. For general $\chi$, one argues with functional calculus using that $t\mapsto \chi(t)-t(1+t^2)^{-1/2}$ is a $C_0$-function. The finite-summability statement can be found in \cite[Proposition 3.14]{CGRS11}. We refer to the process 
$$(\mathcal{A},\mathcal{H},D)\mapsto (\pi,\mathcal{H},D(1+D^2)^{-1/2}),$$
as the bounded transform.

\begin{remark}
The assumption that there exists a $*$-algebra $\mathcal{A}_0$ with $\pi(\mathcal{A}_0)\subseteq \mathcal{A}^2$ is made for technical reasons. Such an algebra exists as soon as $\pi^{-1}(\mathcal{A})\subseteq A$ is dense. The latter condition is not automatic, not even for spectral triples arising from smooth, compact manifolds as shown in \cite[Appendix A]{kkbordism}.
\end{remark}

\begin{definition}
An unbounded cycle of order $m>0$ on a $C^*$-algebra $A$ is a pair $(\mathcal{H},D)$ where $\mathcal{H}$ is a separable Hilbert space carrying a representation $\pi:A\to \mathbb{B}(\mathcal{H})$ and $D$ is a self-adjoint operator on $\mathcal{H}$ fitting into a spectral triple $(\mathcal{A},\mathcal{H},D)$ of order $m$ with $\pi^{-1}(\mathcal{A})\subseteq A$ dense.
\end{definition}

Let us return to the spectral triple $(C^\infty_c(M),L^2(M,E),D)$ of Proposition \ref{adlknadlkna}.

\begin{prop}
\label{adlknadlknawithf}
Let $X$ be a metric space. Assume that $(C^\infty_c(M),L^2(M,E),D)$ is a spectral triple of order $m=1$ as in Proposition \ref{adlknadlkna} and $f:M\to X$ a proper, locally Lipschitz mapping. The representation $f^*:C_0(X)\to \mathbb{B}(L^2(M,E))$ defined from pullback along $f$ defines an unbounded cycle on $C_0(X)$ that we denote by $(_{f^*}L^2(M,E),D)$.

Moreover, if $\Gamma$ is a discrete group acting locally Lipschitz continuously on $X$ and on $M$ as in Proposition \ref{adlknadlkna}, then whenever $f$ is $\Gamma$-equivariant the integrated representation $f^*:C_0(X)\rtimes \Gamma\to \mathbb{B}(L^2(M,E))$ defines an unbounded cycle on $C_0(X)\rtimes \Gamma$ that we by an abuse of notation also denote by $(_{f^*}L^2(M,E),D)$
\end{prop}

\begin{proof}
Whenever $X$ is a metric space and $f:M\to X$ is a proper, locally Lipschitz mapping, $f^*:C_0(X)\to C_0(M)$ satisfies that $\pi(\mathrm{Lip}_c(X))\subseteq \mathrm{Lip}_c(M)$. In light of Remark \ref{extendtolip} we see that $(_{f^*}L^2(M,E),D)$ is an unbounded cycle by the definition applied to the spectral triple $(\mathcal{A},\mathcal{H},D)=(\mathrm{Lip}_c(M),L^2(M,E),D)$, or the spectral triple $(\mathcal{A},\mathcal{H},D)=(\mathrm{Lip}_c(M)\rtimes_{\rm alg}\Gamma,L^2(M,E),D)$ in the equivariant case.
\end{proof}

We now use the bounded cycles defined above in Definition \ref{bddcycle} to define the $K$-homology of a $C^*$-algebra. We first need to separate between the possible parities of cycles.

\begin{definition}
An even bounded (or unbounded) cycle $(\mathcal{H},F)$ (or $(\mathcal{H},D)$) on a $C^*$-algebra $A$ is a bounded (or unbounded) cycle where $\mathcal{H}$ carries the structure of a $\Z/2$-graded Hilbert space in which $D$ is odd and $A$ acts as even operators. 
 
An odd bounded (or unbounded) cycle $(\mathcal{H},F)$ (or $(\mathcal{H},D)$) on a $C^*$-algebra $A$ is a bounded (or unbounded) cycle without further structure.
\end{definition}

\begin{remark}
We see that the constructions and examples above readily can be carried out with a grading. For instance, in Proposition \ref{adlknadlknawithf} $(_{f^*}L^2(M,E),D)$ is an even cycle if $E\to M$ is a $\Z/2$-graded vector bundle in which the differential expression $\slashed{D}$ is odd. In other words, decomposing $E=E_+\oplus E_-$ into its graded components
$$ \slashed{D}=\begin{pmatrix}0& \slashed{D}_+\\\slashed{D}_-& 0\end{pmatrix},$$
where $\slashed{D}_\pm:C^\infty(M,E_\mp)\to C^\infty(M,E_\pm)$ is elliptic with formal adjoint $\slashed{D}_\mp$. Furthermore, Baaj-Julg's theorem \ref{baajjulgthm} respects parity in that it maps even unbounded cycles to even bounded cycles.
\end{remark}

We define an operator homotopy of even/odd bounded cycles on $A$ as a pair $(\mathcal{H},\pmb{F})$ where $\mathcal{H}$ is a Hilbert space with a representation of $A$, and $\pmb{F}:[0,1]\to \mathbb{B}(\mathcal{H})$ is a norm-continuous path such that $(\mathcal{H},\pmb{F}(t))$ is an even/odd bounded cycle on $A$ for any $t\in [0,1]$. If this is the case, we write this operator homotopy as 
$$(\mathcal{H},\pmb{F}(0))\sim_{\rm oh}(\mathcal{H},\pmb{F}(1)).$$
The reader can note that if $(\mathcal{H},F)$ is an even/odd cycle on $A$ and $K$ is a locally compact operator on $\mathcal{H}$ (i.e. $aK,Ka\in \mathbb{K}(\mathcal{H})$), then if the cycle is odd or $K$ is odd when the cycle is even, the expression
$$\pmb{F}(t):=F+tK$$
defines an operator homotopy 
$$(\mathcal{H},F)\sim_{\rm oh}(\mathcal{H},F+K).$$
We say that an even/odd bounded cycle $(\mathcal{H},F)$ on $A$ is degenerate if for any $a\in A$
$$[F,a]=a(F^*-F)=a(F^2-1)=0.$$
Two cycles $(\mathcal{H},F)$ and $(\mathcal{H}',F')$ are said to be unitarily equivalent if there is a unitary $U:\mathcal{H}\to \mathcal{H}'$ with $UFU^*=F'$.

On the set of unitary equivalence classes of even/odd bounded cycles, there is an abelian semigroup structure defined from direct sum
$$(\mathcal{H},F)+(\mathcal{H}',F'):=(\mathcal{H}\oplus \mathcal{H}',F\oplus F').$$ 
On the abelian semigroup of unitary equivalence classes of even/odd bounded cycles, we define the following equivalence relation: two cycles $(\mathcal{H},F)$ and $(\mathcal{H}',F')$ are declared equivalent $(\mathcal{H},F)\sim (\mathcal{H}',F')$ if there are degenerate cycles $(\mathcal{H}_0,F_0)$ and $(\mathcal{H}'_0,F'_0)$ and an operator homotopy 
$$(\mathcal{H},F)+ (\mathcal{H}_0,F_0)\sim_{\rm oh} (\mathcal{H}',F')+ (\mathcal{H}_0',F_0').$$

\begin{definition}
The even $K$-homology $K^0(A)$ of a $C^*$-algebra $A$ is the abelian group of equivalence classes of unitary equivalence classes of even bounded cycles under $\sim$.
 
The odd $K$-homology $K^1(A)$ of a $C^*$-algebra $A$ is the abelian group of equivalence classes of unitary equivalence classes of odd bounded cycles under $\sim$.

We write $K^*(A):=K^0(A)\oplus K^1(A)$ viewed as a $\Z/2$-graded abelian group.
\end{definition}

The $K$-homology group is well studied \cite{Blackadar, higsonroe, knudsenjensen}. One can easily show that it is contravariantly functorial in the $C^*$-algebra. That is, if $\alpha:A\to B$ is a $*$-homomorphism the map defined at the level of cycles 
$$\alpha:(\pi,\mathcal{H},F)\mapsto (\pi\circ \alpha,\mathcal{H},F),$$
induces a $\Z/2$-graded additive map $\alpha^*:K^*(B)\to K^*(A)$. Moreover, $K$-homology is homotopy invariant in the sense that if 
$$\alpha=(\alpha_t)_{t\in [0,1]}:A\to C([0,1],B),$$
is a $*$-homomorphism then 
$$\alpha_0^*=\alpha_1^*:K^*(B)\to K^*(A).$$

We next recall a result of Baaj-Julg showing that the bounded transform studied in Theorem \ref{baajjulgthm} is surjective at the level of cycles modulo locally compact perturbations. For an unbounded cycle $(\mathcal{H},D)$, we write 
$$[(\mathcal{H},D)]:=[(\mathcal{H},D(1+D^2)^{-1/2})],$$
for the $K$-homology class of its bounded transform. We refer to $[(\mathcal{H},D)]$ as  the $K$-homology class of $(\mathcal{H},D)$.

\begin{thm}[Baaj-Julg \cite{baajjulg}]
\label{baajjulgthm2}
Assume that $A$ is a separable $C^*$-algebra. Given a bounded cycle $(\mathcal{H},F)$ on $A$, there is a self-adjoint operator $D$ on $\mathcal{H}$ with the properties
\begin{itemize}
\item $D$ is odd if $(\mathcal{H},F)$ is even;
\item $a(D(1+D^2)-F)\in \mathbb{K}(\mathcal{H})$ for any $a\in A$;
\item $(\mathcal{H},D)$ is an unbounded cycle on $A$.
\end{itemize}
In particular, for any $x\in K^*(A)$ there is an unbounded cycle $(\mathcal{H},D)$ on $A$ such that 
$$x=[(\mathcal{H},D)].$$
\end{thm}

\begin{problem}
\label{summabilityprob}
Let $A$ be a separable $C^*$-algebra. \newline
\noindent {\bf Bounded version}: Given $x\in K^*(A)$, is there a finitely summable bounded cycle $(\mathcal{H},F)$ on $A$ such that 
$$x=[(\mathcal{H},F)]?$$

\noindent{\bf Unbounded version} Given $x\in K^*(A)$, is there a finitely summable unbounded cycle $(\mathcal{H},D)$ on $A$ such that 
$$x=[(\mathcal{H},D)].$$

\end{problem}

Problem \ref{summabilityprob} is motivated both by the fundamental question of when cycles can be represented by objects analogous to finite-dimensional cycles, but also by the computational utility of finite-dimensional cycles. Elements of $K$-homology $K^*(A)$ define invariants on the $K$-theory $K_*(A)$ via the index pairing $K_*(A)\times K^*(A)\to \Z$, see more in \cite{connesncdg,redbook,connesmoscovici,levithesis}. Finitely summable cycles define cyclic cocycles \cite{connesncdg,redbook} explicitly computing the index pairing via an Atiyah-Singer like formula, and in the unbounded setting one can under additional technical assumptions localize said cyclic cocycle to Connes-Moscovici's local index formula \cite{connesmoscovici}.

Let us mention a number of examples where Problem \ref{summabilityprob} has  positive solution. These results are all results of a type showing that certain classes of examples exhaust $K$-homology groups.

\begin{example}
Assume that $X$ is a finite CW-complex. We fix a distance function on $X$ that on each cell is locally bi-Lipschitz equivalent to an open set in Euclidean space. We write $\mathrm{Lip}(X)\subseteq C(X)$ for the dense subalgebra of Lipschitz functions. It holds that any $K$-homology class $x\in K^*(C(X))$ can be represented by a finitely summable unbounded cycle of the form $(_{f^*}L^2(M,E\otimes S_M),D_E)$ as in Proposition \ref{adlknadlknawithf}, where $M$ is a compact spin$^c$–manifold, $f:M\to X$ is a Lipschitz map and $E\to M$ is a complex vector bundle.

This fact follows from Baum-Douglas model for $K$-homology of topological spaces \cite{BaumDouglas1982}, in which cycles are built up from triples $(M,E,f)$ as above. It was proven in \cite{BaumHigsonSchick2007} that under an appropriate equivalence relation on the Baum-Douglas cycles $(M,E,f)$, it forms a generalized homology theory $K_*^{\rm geo}(X)$ on topological spaces. Baum, Higson and Schick proved that, for a finite CW-complex, the analytic assembly map from geometric $K$-homology to analytic $K$-homology is an isomorphism \cite{BaumHigsonSchick2007}. Under this isomorphism the geometric cycle $(M,E,f)$ is sent to the analytic $K$-homology class $[(f^*L^2(M,E\otimes S_M),D_E)]$. Consequently every class in $K^*(C(X))$ admits a finitely summable unbounded representative of this form.
\end{example}

\begin{example}
If $A$ is a purely infinite $C^*$-algebra, there are no non-trivial finitely summable unbounded cycles by the Connes' tracial obstruction \cite{ConnesETDS,GRU}. Nevertheless, there are several results showing that the bounded case of Problem \ref{summabilityprob} has a positive solution any $x\in K^*(A)$ for large classes of purely infinite $C^*$-algebras $A$. We mention Cuntz-Krieger algebras \cite{GM15}, Smale spaces \cite{Ger}, and hyperbolic groups \cite{EN}. 
\end{example}

\section{$KK$-theory}
\label{seckksda}

A fundamental tool for working with and understanding spectral noncommutative geometry is $KK$-theory. The Kasparov product on the $KK$–groups provides topological invariants on $K$-theory and $K$-homology. In its constructive, unbounded, form the Kasparov product provides a technique for manipulating and constructing spectral noncommutative geometries. While $KK$-theory historically has played an important role in the study of the Baum-Connes conjecture and its applications, we here take a converse approach and use ideas from the primordial method to prove the Baum-Connes conjecture as a means to construct well behaved cycles representing $K$-homology classes on crossed products. We first recall basic elements of $KK$-theory \cite{Blackadar,kaspinv,knudsenjensen,thebeastofmesland}, in its bounded and unbounded picture, after which we discuss the important points of the Kasparov product and finally we recall the notion of the $\gamma$-element and examples thereof \cite{kaspinv,tu00}.

We consider two separable $C^*$-algebras $A$ and $B$ with actions of a discrete group $\Gamma$. We often use the notion of a $B$-Hilbert $C^*$-module that extends Hilbert spaces (where $B=\C$) an analytically flexible notion of module to general $C^*$-algebras, see \cite{lancebook} for an introduction. On a $B$-Hilbert $C^*$-module $\mathpzc{E}$ we write $\End_B^*(\mathpzc{E})$ for the $C^*$-algebra of $B$-linear, adjointable operators on $\mathpzc{E}$. We write $\mathbb{K}_B(\mathpzc{E})\subseteq \End_B^*(\mathpzc{E})$ for the closed ideal densely spanned the rank one operators $T_{\xi,\eta}:x\mapsto \xi\langle \eta,x\rangle_{\mathpzc{E}}$, where $\xi$ and $\eta$ ranges over $\mathpzc{E}$. For densely defined operators $D:\mathpzc{E}\dashrightarrow \mathpzc{E}$ one defines notions such as closed operators, adjoint operators and self-adjoint operators similarly to for Hilbert spaces \cite{lancebook}. We often need to assume also that the densely defined operator $D:\mathpzc{E}\dashrightarrow \mathpzc{E}$ is regular in order to extend the theory of densely defined operators on Hilbert spaces, the definition of regularity can be found in \cite{lancebook} or via the so called local/global characterization \cite{leschkaadlocal}.

\begin{definition} 
\label{adljnadjlknad}
Let $A$ and $B$ be two separable $C^*$-algebras with actions of a discrete group $\Gamma$. 
\begin{itemize}
\item A bounded $\Gamma$-equivariant $(A,B)$-cycle is a pair $(\mathpzc{E},F)$ where $\mathpzc{E}$ is a countably generated, $\Gamma$-equivariant $B$-Hilbert $C^*$-module carrying a representation $A\to \End_B^*(\mathpzc{E})$ and $F\in \End_B^*(\mathpzc{E})$ is an operator satisfying for any $a\in A$ and $g\in \Gamma$ that
$$[F,a], \; a(F^*-F), \;a(F^2-1), \; a(g(F)-F)\in \mathbb{K}_B(\mathpzc{E}).$$
\item An unbounded $\Gamma$-equivariant $(A,B)$-cycle of order $m>0$ is a pair $(\mathpzc{E},D)$ where $\mathpzc{E}$ is a countably generated, $\Gamma$-equivariant $B$-Hilbert $C^*$-module carrying a representation $A\to \End_B^*(\mathpzc{E})$ and $D:\mathpzc{E}\dashrightarrow \mathpzc{E}$ is a selfadjoint, regular operator satisfying there is a dense $*$-subalgebra $\mathcal{A}\subseteq A$ such that for any $a\in \mathcal{A}$ and $g\in \Gamma$ we have
\begin{enumerate}
\item $a(i+D)^{-1}\in \mathbb{K}_B(\mathpzc{E})$;
\item $a\Dom(D)\subseteq \Dom(D)$ and $[D,a](1+D^2)^{-1/2+1/2m}$ is norm bounded.
\item $g\Dom(D)\subseteq \Dom(D)$ with $a(g(D)-D)(1+D^2)^{-1/2+1/2m}$ and $(g(D)-D)a(1+D^2)^{-1/2+1/2m}$ extending to bounded operators.
\end{enumerate}
\end{itemize}
\end{definition}

Similarly to the previous subsection, we can define unitary equivalence of bounded and unbounded cycles. We can also define direct sums of bounded cycles, and direct sum defines a binary operation also on the unbounded cycles when fixing the subalgebra $\mathcal{A}$. The set of unitary equivalence classes of (un-) bounded $(A,B)$-cycles depends covariantly on $B$ and contravariantly on $A$ (where the functoriality of unbounded cycles require that the dense subalgebras are preserved).

A homotopy of bounded $\Gamma$-equivariant $(A,B)$-cycle is a $\Gamma$-equivariant $(A,C([0,1],B))$-cycle $(\mathpzc{E},F)$. Writing $\mathrm{ev}_t:C([0,1],B)\to B$ for the evaluation mappings, we declare the endpoints $(\mathrm{ev}_0)_*(\mathpzc{E},F)$ and $(\mathrm{ev}_1)_*(\mathpzc{E},F)$ as homotopy equivalent bounded $\Gamma$-equivariant $(A,B)$-cycles. Similarly, following \cite{kaadhomoun}, we can also define a homotopy of unbounded $\Gamma$-equivariant $(A,B)$-cycles as a $\Gamma$-equivariant unbounded $(A,C([0,1],B))$-cycle $(\mathpzc{E},D)$, and declare the endpoints $(\mathrm{ev}_0)_*(\mathpzc{E},D)$ and $(\mathrm{ev}_1)_*(\mathpzc{E},D)$ as homotopy equivalent unbounded $\Gamma$-equivariant $(A,B)$-cycles. As in the previous subsection, we also define even $\Gamma$-equivariant $(A,B)$-cycles or homotopies as those where $\mathpzc{E}$ is $\Z/2$-graded, $A$ and $\Gamma$ acts evenly and $F$ (or $D$ in the unbounded case) is odd. Homotopy defines an equivalence relation on the set of unitary equivalence classes of odd or even $\Gamma$-equivariant $(A,B)$-cycles, leading us to the definition of the $KK$-groups.

\begin{definition}
Let $A$ and $B$ be two separable $C^*$-algebras with actions of a discrete group $\Gamma$. The even $KK$-group $KK_0^\Gamma(A,B)$ is the abelian group of homotopy classes of unitary equivalence classes of even bounded $\Gamma$-equivariant $(A,B)$-cycles.
 
The odd $KK$-group $KK_1^\Gamma(A,B)$ is the abelian group of homotopy classes of unitary equivalence classes of odd bounded $\Gamma$-equivariant $(A,B)$-cycles

We write $KK_*^\Gamma(A,B):=KK_0^\Gamma(A,B)\oplus KK_1^\Gamma(A,B)$ viewed as a $\Z/2$-graded abelian group. If $\Gamma$ is the trivial group, we drop it from the notation.
\end{definition}

It is well known that 
$$KK_*^\Gamma(A,\C)\cong KK_*(A\rtimes\Gamma,\C)=K^*(A\rtimes\Gamma), \quad\mbox{and in particular} \quad KK_*(A,\C)=K^*(A).$$
The first isomorphism is implemented at the level of cycles by $(\mathcal{H},F)\mapsto (\mathcal{H},F)$ where the left hand side is viewed as a $\Gamma$-equivariant $A$-cycle and the right-hand side is viewed as an $A\rtimes\Gamma$-cycle via the integrated representation. The second isomorphism reduces to showing $KK_*(A,\C)=K^*(A)$ which is proven \cite{higsonroe} by showing that the homotopy relation coincides with the equivalence relation generated by operator homotopy and degenerate equivalence. We moreover have the isomorphism \cite[Chapter 17.5]{Blackadar}
$$KK_*(\C,B)\cong K_*(B),$$
via an index mapping.

\begin{remark}
In fact, for $\Gamma$ being the trivial group, both theorems of Baaj-Julg above, Theorems \ref{baajjulgthm} and \ref{baajjulgthm2}, extend to the bivariant setting. That is, any unbounded $(A,B)$-cycle $(\mathpzc{E},D)$ defines a class 
$$[(\mathpzc{E},D)]:=[(\mathpzc{E},D(1+D^2)^{-1/2})]\in KK_*(A,B),$$
and any class $x\in KK_*(A,B)$ can be represented by an unbounded $(A,B)$-cycle. Kaad \cite{kaadhomoun} goes further and shows that if we fix the smooth structure $\mathcal{A}\subseteq A$ as a countable $*$-subalgebra, then $KK_*(A,B)$ is isomorphic to the set of homotopy classes of unbounded $(A,B)$-cycles with $\mathcal{A}$ fixed. Similar results were also obtained independently by van den Dungen-Mesland \cite{meslsdlda}.
\end{remark}

We have a construction of descent of $\Gamma$-equivariant $(A,B)$-cycles into $(A\rtimes \Gamma,B\rtimes \Gamma)$-cycles
$$(\mathpzc{E},F)\mapsto (\mathpzc{E}\rtimes \Gamma,F\rtimes \Gamma)\quad\mbox{and}\quad(\mathpzc{E},D)\mapsto (\mathpzc{E}\rtimes \Gamma,D\rtimes \Gamma).$$
This produces a descent mapping 
$$KK_*^\Gamma(A,B)\to KK_*(A\rtimes \Gamma,B\rtimes \Gamma).$$

One of the most important features of $KK$-theory is the Kasparov product. Kasparov \cite{kaspinv} proved that for any three $\Gamma-C^*$-algebras $A$, $B$, and $C$ there is a product operation 
$$\otimes_B:KK_*^\Gamma(A,B)\times KK_*^\Gamma(B,C)\to KK_*^\Gamma(A,C),$$
which is associative and functorial. We shall neither define nor construct the Kasparov product, but only use it as a black box. What we will need is a method of computing it in the unbounded picture. This approach was after initial results of Kucerovsky, pioneered by Mesland \cite{thebeastofmesland} and later further developed \cite{leschkaadsf,leschmes} in the community. We recall the following generalization \cite[Theorem 6.3.4]{thebeastofmesland}. It follows from \cite[Theorem A.7]{GM15} (extending Kucerovsky's theorem \cite{dantheman} to higher order) and Kaad-Lesch's result on sums of self-adjoint operators \cite{leschkaadlocal,leschkaadsf}.

\begin{thm}
\label{alkdnajdn}
Let $A$ and $B$ be two separable $C^*$-algebras, $(\mathpzc{E},S)$ an unbounded $(A,B)$-cycle of order $m\geq 1$ and $(\mathcal{B},\mathcal{H},D)$ a spectral triple of order $1$ with $\mathcal{B}\subseteq B$ dense and the inclusion   $\mathcal{B}\hookrightarrow B$ continuous in norm. Assume that $\mathcal{E}\subseteq  \mathpzc{E}$ is a dense subspace equipped with a linear map
$$\nabla:\mathcal{E}\to \mathcal{E}\otimes^{\rm alg}_\mathcal{B}\mathbb{B}(\mathcal{H}),$$
such that 
\begin{enumerate}
\item $\mathcal{E}$ is a $\mathcal{B}$-module and $\nabla$ is a connection in the sense that for $\xi\in \mathcal{E}$, $b\in \mathcal{B}$ we have 
$$\nabla(\xi b)=\nabla(\xi)b+\xi\otimes [D,b].$$
\item The densely defined operator on $\mathpzc{E}\otimes_B \mathcal{H}$ given by
$$\xi\otimes x\mapsto \nabla(\xi)x+\xi\otimes Dx, \quad \xi\in \mathcal{E}, \; x\in \Dom(D),$$
is essentially self-adjoint and its closure $1\otimes_\nabla D$ satisfies 
\begin{itemize}
\item there is a dense $*$-subalgebra $\mathcal{A}$ such that for $a\in \mathcal{A}$, 
\begin{align*}
a\Dom(1\otimes_\nabla D)&\subseteq \Dom(1\otimes_\nabla D),\quad\mbox{and}\\
& [1\otimes_\nabla D,a]\quad\mbox{is bounded, and}
\end{align*}
\item there is a core $\tilde{E}$ for the operator $1\otimes_\nabla D$ such that $(S\otimes 1\pm i)^{-1}\tilde{E}\subseteq \Dom ((S\otimes 1)(1\otimes_\nabla D))\cap \Dom(1\otimes_\nabla D)\cap \Dom(S\otimes 1)$ and the anti-commutator $[1\otimes_\nabla D,S\otimes 1]_+$ satisfies that 
$$[1\otimes_\nabla D,S\otimes 1]_+(S\otimes 1\pm i)^{-1}:\tilde{E}\to \mathpzc{E}\otimes_B \mathcal{H},$$
is norm bounded.
\end{itemize}
\end{enumerate}
Then it holds that $(\mathpzc{E}\otimes_B \mathcal{H},S\otimes 1+1\otimes_\nabla D)$ is an unbounded Fredholm module of order $m$ on $A$ and 
$$[(\mathpzc{E}\otimes_B \mathcal{H},S\otimes 1+1\otimes_\nabla D)]=[(\mathpzc{E},S)]\otimes_B[(\mathcal{B},\mathcal{H},D)].$$
\end{thm}

\section{$\gamma$-elements}
\label{subsecgamma}

We now proceed with a pedestrian overview of Kasparov's idea of using $\gamma$-elements to prove the Baum-Connes conjecture, where we do not explicitly use the Baum-Connes conjecture but rather the $\gamma$-element itself. The $\gamma$-element of a discrete group $\Gamma$ is, when it exists, an element $\gamma\in KK^\Gamma(\C,\C)$ that in a sense regularizes $K$-theoretical behaviour in that it in an appropriate categorical sense projects onto the part coming from proper actions \cite{meyernest}. The definition of a $\gamma$-element was recalled above at the beginning of Section \ref{subsecintro}, see also \cite{kaspinv,nishisplit,tu00}. The approach we take to $\gamma$-elements relies on Kasparov's notion of special $\Gamma$-manifolds \cite{kaspinv} that we recall next.

Following \cite{kaspinv}, a complete, second countable, Riemannian manifold $M$ with a smooth, isometric action of $\Gamma$ is called a $\Gamma$-manifold. We will tacitly assume that $M$ is connected. We note that if $\Gamma$ acts properly, cocompactly and smoothly on a second countable manifold $M$, there is always a complete Riemannian metric in which the $\Gamma$-action is isometric. We us the term proper $\Gamma$-manifold for a $\Gamma$-manifold where the action is proper. We write $C_0(M,\C\ell(M))$ for the $\Z/2$-graded $C^*$–algebra of sections of the bundle of complex $\Z/2$-graded Clifford algebras of the cotangent bundle on $M$. The isometric action of $\Gamma$ induces a $\Gamma$-action on the $C^*$-algebra $C_0(M,\C\ell(M))$. Since $M$ is complete, the Hodge-Dirac operator $\rd+\rd^*$ acting on the exterior algebra bundle of $M$ is essentially self-adjoint (cf. Remark \ref{chernoffrem}) with self-adjoint closure $D$ and the graded action of $C_0(M,\C\ell(M))$ on $L^2(M,\wedge^*T^*M\otimes \C)$ makes $(L^2(M,\wedge^*T^*M\otimes \C), D)$ into a $\Gamma$-equivariant, graded unbounded cycle for  $C_0(M,\C\ell(M))$, defining a class in $K^0_\Gamma(C_0(M,\C\ell(M)))=KK_0^\Gamma(C_0(M,\C\ell(M)),\C)$. Following \cite{kaspinv}, we denote this class by $[\rd_M]$.

\begin{definition}\label{special}
A $\Gamma$-manifold $M$ is called special if there is an element 
$$\eta_M\in KK_0^\Gamma(\C,C_0(M,\C\ell(M))),$$ 
such that 
$$[\rd_M]\otimes_\C\eta_M=1\in KK_0^\Gamma(C_0(M,\C\ell(M))),C_0(M,\C\ell(M))).$$
\end{definition}

\begin{remark}
\label{rempos}
It was proven in \cite[Page 186]{kaspinv} that a $\Gamma$-manifold with non-positive sectional curvature produces a special $\Gamma$-manifold. A fact we shall use is that on  \cite[Page 186]{kaspinv}, it is in fact proven that $\eta$ can be represented by a particularly nice unbounded $(\C,C_0(M,\C\ell(M)))$-cycle. We write $\rho$ for the geodesic distance function on $M$. We fix a point $x_0\in M$. Define the unbounded multiplier $\delta\in C^\infty(M,\C\ell(M))$ of $C_0(M,\C\ell(M))$ as 
$$\delta(x):=\rho(x_0,x)\rd_x\rho(x_0,x).$$
We have that $\delta(1+\delta^2)^{-1/2}=\rd_x(1+\rho(x_0,x)^2)^{1/2}$ coincides with the bounded multiplier $\xi_{x_0}$ appearing on \cite[Page 186]{kaspinv}. Therefore, $(C_0(M,\C\ell(M)),\delta)$ is a graded, $\Gamma$-equivariant, unbounded $(\C,C_0(M,\C\ell(M)))$-cycle of order $m=1$ with 
$$\eta_M:=[((C_0(M,\C\ell(M)),\delta)]\in KK_0^\Gamma(\C,C_0(M,\C\ell(M))),$$
fulfilling the property in the definition of a special $\Gamma$-manifold.
\end{remark}

\begin{prop}
\label{betarep}
Assume that $M$ is a proper, special $\Gamma$-manifold then with $\beta:=\eta_M \in KK_0^\Gamma(\C,C_0(M,\C\ell(M)))$ and $\alpha:=[\rd_M]\in  KK_0^\Gamma(C_0(M,\C\ell(M)),\C)$, the element $\gamma:=\beta \otimes_{C_0(M,\C\ell(M))}\alpha\in KK^\Gamma_0(\C,\C)$ defines a $\gamma$-element for $\Gamma$, which is geometric in the sense of Definition \ref{geomegamadma}.

If moreover $M$ is a proper $\Gamma$-manifold with non-positive sectional curvature, the $\gamma$-element constructed above has a geometric $\beta$-factor of order $m=1$ defined from the graded, $\Gamma$-equivariant, unbounded $(\C,C_0(M,\C\ell(M)))$-cycle $(C_0(M,\C\ell(M)),\delta)$. 
\end{prop}

\begin{example}
\label{aodknaldn}
Let $\Gamma$ be a finitely generated group of polynomial growth. By Gromov's theorem, $\Gamma$ is virtually nilpotent \cite{Gromov1981}. A virtually nilpotent group $\Gamma$ admits a proper cocompact affine model on a simply connected nilpotent Lie group $\mathsf{N}$; see, for example, \cite{Dekimpe2003}. By \cite[Theorem 5.7]{kaspinv}, we have that $\mathsf{N}$ is a special $\Gamma$-manifold. However, in general $\mathsf{N}$ does not have non-positive sectional curvature so the construction in Remark \ref{rempos} is not sufficient to allow us to deduce that $\gamma_\Gamma$ has a geometric $\beta$-factor. For the proof of Theorem \ref{mainthm3} we shall need the following strengthening of the dual-Dirac construction.
\end{example}

\begin{lemma}
\label{lemaodknaldn}
Let $\Gamma$ be a finitely generated nilpotent group. Then $\gamma_\Gamma$ admits a geometric $\beta$-factor. Indeed, $\beta$ is represented by an unbounded $\Gamma$-equivariant cycle of order $m\geq 1$ of the form $(C_0(\mathsf{N},\mathcal{E}),V)$ where $\mathsf{N}$ is a nilpotent Lie group containing a finite quotient of $\Gamma$ as a lattice, $\mathcal{E}\to \mathsf{N}$ is a hermitean vector bundle and $V$ is a polynomial endomorphism valued section with $V^2$ scalar valued and $|\nabla V|\lesssim (1+|V|)^{\frac{m-1}{m}}$ for $m>0$ being the graded order of the polynomial $V$. 
\end{lemma}

\begin{proof}
By \cite[Theorem 5.2.7]{crusoe}, $\Gamma$ has a finite normal subgroup $F$ such that $\Gamma/F$ is torsion-free. We can in particular assume that $\Gamma$ is a torsion-free, finitely generated nilpotent group. We take $\mathsf{N}$ as the Malcev completion of $\Gamma$, so the uniquely determined simply connected nilpotent Lie group containing $\Gamma$ as a lattice. It suffices to show that the $\beta$-element of $\mathsf{N}$ has the desired form. 

Following the proof of \cite[Theorem 5.7]{kaspinv}, we can construct $\beta$ inductively as a Kasparov product using the lower central series of $\mathsf{N}$. This is carried out in \cite[Section 7]{FGM} showing that $\beta$ is represented by a cycle of the form $(C_0(\mathsf{N},\mathcal{E}),V)$, for $\mathcal{E}\to \mathsf{N}$ a hermitean vector bundle and $V$ a polynomial endomorphism valued section with $V^2$ scalar valued. Indeed, $V^2$ is a power of a quasi-norm and polynomial homogeneous of even order $2m$. Since $V^2$ is a quasi-norm, $|\nabla V|\lesssim (1+|V|)^{\frac{m-1}{m}}$ follows.
\end{proof}

\section{Finite-summability and Weyl laws}
\label{alkdnaldn}

Above in Proposition \ref{adlknadlkna}, we saw how a first-order elliptic operator $\slashed{D}$ on a complete manifold $M$ gives rise to a finitely summable spectral triple on the non-unital algebra $C^\infty_c(M)$. Our results in this section lead to results proving that under appropriate conditions on a potential $V$ and assuming $M$ has bounded geometry, the operator $\slashed{D}+V$ produces a finitely summable spectral triple on bounded Lipschitz functions on $M$; a somewhat different type of result as in this case the finite-summability condition is global. More precisely, we review existing work in spectral theory that we use to prove that for a Dirac operator $\slashed{D}$ on a bounded geometry and a self-adjoint locally Lipschitz potential $V$, satisfying appropriate conditions on $\nabla V$ and anticommuting with the symbol of $\slashed{D}$, the operator $\slashed{D}+V$ has its resolvent in a Schatten class $\mathcal{L}^p$ (for some $p$) if and only if $(i+V)^{-1}$ belongs to some Lebesgue class $L^q$ (for some $q$).

\begin{definition}
A complete Riemannian manifold $M$ with metric $g$ is said to have bounded geometry if it has positive injectivity radius and its curvature and all its covariant derivatives are bounded.

A hermitean vector bundle with a hermitean connection $(E,\nabla_E)$ on a manifold $M$, of bounded geometry, is said to have bounded geometry if its curvature and all its covariant derivatives are bounded. 
\end{definition}

\begin{example}
If $M$ is a complete Riemannian manifold with a proper, cocompact and isometric action of a discrete group $\Gamma$ then $M$ has bounded geometry by the following argument. Because the action is isometric, the injectivity-radius function on $M$ is $\Gamma$-invariant. It is positive and continuous on $M$, and therefore has
a positive lower bound on a compact set $F\subset M$ with $\Gamma F=M$. The curvature tensor and all of its covariant derivatives are likewise $\Gamma$-invariant, so their norms are bounded on $F$ and hence on all of $M$. Thus $M$ has bounded geometry.

Similarly, if $E\to M$ is an equivariant hermitean vector bundle with an equivariant hermitean connection $(E,\nabla_E)$ on the manifold $M$ with a proper, cocompact and isometric action of a discrete group $\Gamma$, then $(E,\nabla_E)$ has bounded geometry.
\end{example}

We note that if $M$ and $E\to M$ have bounded geometry, there is a canonically associated scale of Sobolev spaces $H^s(M,E)$ and any differential operator $D$ of bounded geometry and of order $m$ acts continuously $D:H^s(M,E)\to H^{s-m}(M,E)$. Moreover, any uniformly elliptic operator $D$ of order $m$ extends uniquely from $C^\infty_c(M,E)$ to a closed operator on $L^2(M,E)$ and its domain is $H^m(M,E)$. In particular, the domain of Dirac operators as in Example \ref{akandaojnd} is $H^1(M,S_M\otimes E)$. For more details, see \cite{boundedgeome}

\begin{thm}
\label{finsumthm}
Assume that $M$ is a complete Riemannian spin$^c$ manifold of bounded geometry, with a hermitean vector bundle with a hermitean connection $(E,\nabla_E)$ of bounded geometry. Let $\slashed{D}_E$ denote a spin$^c$-Dirac operator on $M$ twisted by $(E,\nabla_E)$ and write $D_E$ for its self-adjoint closure on $L^2(M,S_M\otimes E)$ (cf. Remark \ref{akandaojnd}). Assume that $V\in C(M,\mathrm{End}(S_M\otimes E))$ is self-adjoint and satisfies
\begin{enumerate}
\item $V$ is locally Lipschitz with $V^2$ scalar valued and for some $\alpha>0$ satisfies $|\nabla V|\lesssim (1+|V|)^{1-\frac{1}{\alpha}}$;
\item $V$ anticommutes with the principal symbol of $\slashed{D}_E$.
\end{enumerate}
Then $D_E+V$ is self-adjoint. Moreover there is a $p$ with $(i+D_E+V)^{-1}\in \mathcal{L}^p(L^2(M,S_M\otimes E))$ iff there is a $q$ with $(i+V)^{-1}\in L^q(M,\mathrm{End}(S_M\otimes E))$.
\end{thm}

We note that the assumption that $V$ is locally Lipschitz with $V^2$ scalar valued and for some $\alpha>0$ satisfies $|\nabla V|\lesssim (1+|V|)^{1-\frac{1}{\alpha}}$ ensures that  $D_E+V$ is self-adjoint by \cite[Theorem 1.1]{leschmes}, see also \cite[Section 7]{leschkaadlocal}. The novel part of Theorem \ref{finsumthm} concerns finite-summability.

\begin{lemma}
\label{reductoscalar}
Let $\alpha>0$. Under the assumptions of Theorem \ref{finsumthm}, the equivalence therein holds as soon as for the Laplace-Beltrami operator $\Delta$ on $M$, the following two statements are equivalent for any positive, locally Lipschitz function $v$ with $|\nabla v|\lesssim |v|^{1-\frac{1}{2\alpha}}$ on $M$:
\begin{enumerate}
\item There is a $p$ with $(1+\Delta+v)^{-1}\in \mathcal{L}^p(L^2(M))$.
\item There is a $q$ with $(1+v)^{-1}\in L^q(M)$.
\end{enumerate}
\end{lemma}

\begin{proof}
Take $V$ as in the assumptions of Theorem \ref{finsumthm} and consider the scalar valued function $v:=V^2$. Then $v$ is locally Lipschitz on $M$ with $|\nabla v|\lesssim |v|^{1-1/2\alpha}$. Clearly, $(1+v)^{-1}$ is $L^q$ for some $q$ if and only if $(i+V)^{-1}\in L^q(M,\mathrm{End}(S_M\otimes E))$ for some $q$. By the Lichnerowicz–Weitzenböck formula \cite{BGV} and the bounded geometry assumption ensures that the two closed forms
$$\mathfrak{q}_{D,V}(U):=\int_M (|\slashed{D}_E U|^2+\langle VU,VU\rangle_{S_M\otimes E)})\rd x,$$
and 
$$\mathfrak{q}_V(U):=\int_M (|\nabla_E U|^2+\langle VU,VU\rangle_{S_M\otimes E)})\rd x,$$
have the same form domain $\{u\in H^1(M,S_M\otimes E): VU\in L^2\}$. Here $\nabla_E$ denotes the twisted spin$^c$-connection on $S_M\otimes E$. Since $V$ anticommutes with the principal symbol of $D_E$, the anticommutator $[D_E,V]_+$ is a zeroth-order endomorphism and
\[\|[D_E,V]_+(x)\|\le C|\nabla V(x)| \le C(1+|V(x)|)^{1-1/\alpha}.\]
For every $\varepsilon>0$ there is $C_\varepsilon>0$ such that $(1+t)^{1-1/\alpha}\le \varepsilon t^2+C_\varepsilon, \qquad t\ge0,$ and consequently,
\[|\langle[D_E,V]_+u,u\rangle| \le \varepsilon\|Vu\|^2+C_\varepsilon\|u\|^2.\]
Using the Lichnerowicz--Weitzenbock formula  \cite{BGV}, after increasing the additive constant if necessary, we obtain for every $\varepsilon>0$ a constant $C>0$ such that for $u$ in the form domain
$$(1-\varepsilon)\mathfrak{q}_V(U)-C\|U\|_{L^2}^2\leq \mathfrak{q}_{D,V}(U)\leq (1+\varepsilon)\mathfrak{q}_V(U)+C\|U\|_{L^2}^2.$$
The proof of the Lemma is therefore complete if there are constants $C_1,C_1',C_2,C_2'>0$ and an additional multiplicity $N\in \N$ such that we have singular value estimates
\begin{equation}
\label{glasest}
C_1\mu_n(\Delta+v)-C_1'\leq \mu_{Nn}(\Delta_E+V^2)\leq C_2\mu_n(\Delta+v)+C_2',
\end{equation}
where $\Delta_E:=\nabla_E^*\nabla_E$ is the Bochner laplacian so $\Delta_E+V^2$ is the operator associated with $\mathfrak{q}_V$. We shall conclude the estimate \eqref{glasest} from quadratic form arguments and the Glazman lemma \cite[Theorem 1.25]{franklapweid}, also known as the min-max principle. The quadratic form estimates will come from a partition of unity argument.

We introduce the closed, positive quadratic form 
$$\mathfrak{q}(u):=\int_M (|\nabla u|^2+v|u|^2)\rd x,$$
defined on the form domain $\Dom(\mathfrak{q}):=\{u\in H^1(M): v|u|^2\in L^1\}$.

Since $M$ and $E$ have bounded geometry, we can for $r $ smaller than the injectivity radius choose the following data
\begin{itemize}
\item A countable cover $(U_j)_{j=1,2,3,\ldots}$ with base points $x_j\in U_j$ and diffeomorphisms defined from the exponential map at $x_j$:
$$\phi_j: B(0,r)\to U_j,$$
where $B(0,r)$ denotes the Euclidean ball of radius $r$, such that for some constant $C$ we can bound 
\begin{equation}
\label{lknlknad}
\|\phi_j\|_{C^1}\leq C,\quad \forall j.
\end{equation}
\item Unitary trivializations $\rho_j: (S_M\otimes E)|_{U_j}\xrightarrow{\sim} B(0,r)\times S_n\otimes \C^r$ for $j=1,2,3,\ldots$ covering $\phi_j$, where $S_n$ denotes the spinors on $\R^n$ and $r$ denotes the rank of $E$, such that for some constant $\tilde{C}$ we can bound 
\begin{equation}
\label{lknlknad2}
\|\rho_j\|_{C^1}\leq \tilde{C},\quad \forall j.
\end{equation}
\item A partition of unity $(\chi_j)_{j=1,2,3,\ldots}$ subordinate to $(U_j)_{j=1,2,3,\ldots}$ such that for some constant $\hat{C}$ we can bound  
\begin{equation}
\label{lknlknad3}
\|\chi_j\|_{C^1}\leq \hat{C},\quad \forall j.
\end{equation}
\end{itemize}
Let $\mu$ denote the Riemannian volume density on $M$, write $\rd m$ for the Lebesgue measure in $\R^n$ and set 
$$c_j:=\sqrt{\frac{\rd((\phi_j^{-1})_*\mu}{\rd m}}.$$
We form the isometries 
$$\mathfrak{v}:L^2(M)\to \bigoplus_{j=1}^\infty L^2(B(0,r)),\quad \mathfrak{v}f:=(\phi_j^*(\chi_j f))_{j=1,2,\ldots}$$
and 
$$\mathfrak{v}_S:L^2(M,S_M\otimes E)\to \bigoplus_{j=1}^\infty L^2(B(0,r),S_n\otimes \C^r),\quad \mathfrak{v}_Sf:=(\rho_j[\phi_j^*(\chi_j f)])_{j=1,2,\ldots}.$$
The reader should note that these isometries on $L^2$ induces continuous embedding with closed ranges
$$\mathfrak{v}:H^s(M)\to \bigoplus_{j=1}^\infty H^s_0(B(0,r)),\quad\mbox{and}\quad \mathfrak{v}_S:H^s(M,S_M\otimes E)\to \bigoplus_{j=1}^\infty H^s_0(B(0,r),S_n\otimes \C^r).$$
We choose a basis $S_n\otimes \C^r\cong \C^N$ and consider the unitary operator 
$$T:=(\mathfrak{v}^{\oplus N})^* \mathfrak{v}_S:L^2(M,S_M\otimes E)\to L^2(M)^N.$$
We note that $T$ induces an invertible operator $H^1(M,S_M\otimes E)\to H^1(M)^N$ and that $T(\Dom(\mathfrak{q}_{V}))\subseteq \Dom(\mathfrak{q}_{v})^{\oplus N}$ and $T^{-1}(\Dom(\mathfrak{q}_{v})^{\oplus N})\subseteq \Dom(\mathfrak{q}_{V})$. A short computation using the bounds \eqref{lknlknad}, \eqref{lknlknad2} and \eqref{lknlknad3} shows that for some $C_1,C_1',C_2,C_2'>0$ we have form estimates for 
$$C_1\mathfrak{q}^{\oplus N}_{v}(Tu_1)-C_1'\|u_1\|^2_{L^2}\leq \mathfrak{q}_{V}(u_1) \quad \mbox{and}\quad \mathfrak{q}_{V}(T^{-1}u_2)\leq C_2\mathfrak{q}^{\oplus N}_{v}(u_2)+C_2'\|u_2\|^2_{L^2},$$
for $u_1\in \Dom(\mathfrak{q}_{V})$ and $u_2\in \Dom(\mathfrak{q}_{v})^{\oplus N}$. From Glazman's lemma, we conclude that the bound \eqref{glasest} holds.
\end{proof}

Theorem \ref{finsumthm} now follows using Lemma \ref{reductoscalar} and the next result.

\begin{lemma}
\label{schattenbounds}
Let $M$ be a manifold of bounded geometry, with $\Delta$ denoting the Laplace-Beltrami operator on $M$. Assume that $\alpha>0$. Let $v$ be a positive, locally Lipschitz function with $|\nabla v|\lesssim (1+|v|)^{1-1/2\alpha}$ on $M$. Then the following conditions are equivalent:
\begin{enumerate}
\item There exists $p\geq 1$ such that
\[(1+\Delta+v)^{-1}\in \mathcal L^p(L^2(M)).\]
\item There exists $q\geq 1$ such that
\[(1+v)^{-1}\in L^q(M). \]
\end{enumerate}
\end{lemma}

\begin{proof}
We can without restriction assume that $v\geq 1$. We first prove that (2) implies (1). Let $k_t(x,y)$ denote the heat kernel of $\Delta$. On a manifold of bounded geometry there is a constant $C_0>0$ such that 
\[
k_t(x,x)\leq C_0t^{-n/2},
\qquad 0<t\leq 1,
\]
see \cite{chegrotay}. The well-known Golden--Thompson inequality therefore gives
\[\Tr(e^{-t(\Delta+v)}) \leq \int_M k_t(x,x)e^{-tv(x)}\,dx.\]
If $(1+v)^{-1}\in L^q(M)$, then, for $0<t\leq 1$, we have $e^{-tv(x)}\lesssim t^{-q}(1+v(x))^{-q}$ and consequently,
\[\Tr(e^{-t(\Delta+v)})\lesssim t^{-n/2-q}, \qquad 0<t\leq 1.\]
Using the integral formula
\[(1+\Delta+v)^{-p}=\frac{1}{\Gamma(p)}\int_0^\infty t^{p-1}e^{-t(1+\Delta+v)}\,dt,\]
we conclude that
\[
\operatorname{Tr}\bigl((1+\Delta+v)^{-p}\bigr)<\infty
\]
for every $p>n/2+q$. This proves (1) assuming (2).

Conversely, suppose that $(1+\Delta+v)^{-1}\in\mathcal L^p(L^2(M))$ for some $p$. If $N(\lambda)$ denotes the eigenvalue counting
function of $1+\Delta+v$, then $N(\lambda)\leq C(1+\lambda)^p$ if $(1+\Delta+v)^{-1}\in\mathcal L^p(L^2(M))$. For $\lambda\geq 1$, we introduce the notation $M_\lambda=\{x\in M:v(x)\leq\lambda\}$ and $w=v^{1/2\alpha}$. Since
\[|\nabla  w| =\frac{|\nabla v|}{ 2\alpha v^{1-1/2\alpha}},\]
the function $w$ is uniformly Lipschitz. Fix $r=c\lambda^{-1/2}$, where $c>0$ is sufficiently small compared with the injectivity radius of $M$. Choose a maximal $2r$-separated family $x_1,\ldots,x_m$ in $M_\lambda$. Bounded geometry implies uniform estimates
\[ c_1r^n\leq \operatorname{vol}(B(x,r)) \leq c_2r^n\]
for sufficiently small $r$, and maximality gives $m\geq c_3\lambda^{n/2}\operatorname{vol}(M_\lambda)$. For each $j$, choose a normalized smooth function $u_j$ supported in $B(x_j,r)$ such that
\[\int_M|\nabla u_j|^2\,dx\leq C_1r^{-2}.\]
The supports may be chosen disjoint. Since $w$ is uniformly Lipschitz and $v(x_j)\leq\lambda$, one has $v(x)\leq C_2\lambda$ on $B(x_j,r)$. Hence we have the bound
\[ \left\langle (1+\Delta+v)u_j,u_j \right\rangle \leq C_3\lambda.\]
The min--max principle \cite{franklapweid} and the established  bound $m\geq c_3\lambda^{n/2}\operatorname{vol}(M_\lambda)$ now implies 
\[ N(C_3\lambda) \geq c_3\lambda^{n/2}\operatorname{vol}(M_\lambda).\]
Since $N(\lambda)\leq C(1+\lambda)^p$, we have that
\[\operatorname{vol}(M_\lambda) \lesssim\lambda^{p-n/2}.\]
We conclude that $(1+v)^{-1}\in L^q(M)$ or every $q>p-n/2$. Therefore (2) follows from (1).
\end{proof}

\begin{cor}
\label{lakdjnadln}
Assume that $M$ is a complete Riemannian spin$^c$ manifold with a proper, cocompact and isometric action of a discrete group $\Gamma$. Assume that $(E,\nabla_E)$ is an equivariant, hermitean vector bundle with a hermitean connection. Let $\slashed{D}_E$ denote a spin$^c$-Dirac operator on $M$ twisted by $(E,\nabla_E)$ and write $D_E$ for its self-adjoint closure on $L^2(M,S_M\otimes E)$ (cf. Remark \ref{akandaojnd}). Assume that $V\in C(M,\mathrm{End}(S_M\otimes E))$ is self-adjoint with $V^2$ scalar valued and satisfies
\begin{enumerate}
\item $V$ for some $m>0$ satisfies that for all $g\in \Gamma$, there is a $C>0$ with 
$$|g^*V-V|\leq C(1+|V|)^{1-1/m};$$
\item $V$ is locally Lipschitz and for some $\alpha>0$ satisfies $|\nabla V|\lesssim (1+|V|)^{1-1/\alpha}$;
\item $V$ anticommutes with the principal symbol of $\slashed{D}_E$.
\end{enumerate}
Then $D_E+V$ is self-adjoint and $\Gamma$-invariant up to perturbations that are relatively bounded by $(1+(D_E+V)^2)^{1/2-1/2m}$. Moreover, there is a $p$ with $(i+D_E+V)^{-1}\in \mathcal{L}^p(L^2(M,S_M\otimes E))$ then $M$ has polynomial growth.
\end{cor}

\begin{proof}
By Theorem \ref{finsumthm}, the assumption $(i+D_E+V)^{-1}\in\mathcal{L}^p$ for some $p$ implies that $(i+V)^{-1}\in L^q(M,\End(S_M\otimes E))$ for some $q$. We define
\[r(x):=(1+\left\|V(x)\right\|_{\End(S_M\otimes E)})^{1/m}. \]
By concavity of $t\mapsto t^{1/m}$ and assumption (1) above, we have for any $g\in \Gamma$ that there is a $C>0$ with $|r(g x)- r(x)|\leq C$ for all $x$. Since the action of $\Gamma$ on $M$ is proper and cocompact, $\Gamma$ is finitely generated; this follows from the Švarc--Milnor lemma \cite{BridsonHaefliger1999}. Fix a finite symmetric generating set $S$ and let $|\gamma|$ denote the associated word length.

For every $s\in S$, the assumption that $V$ is $\Gamma$-invariant up to bounded perturbations implies that there is a constant $C_s$ such that $|r(sx)-r(x)|\leq C_s$, for any $x\in M$. Since $S$ is finite, there is a constant $C>0$ such that
\[r(g x)\leq r(x)+C|g|, \qquad g\in\Gamma,\quad x\in M. \]

Choose a compact set $F\subseteq M$ with $\Gamma F=M$. The family $(g F)_{g\in\Gamma}$ has uniformly bounded multiplicity because the action is proper. Since $r$ is bounded on $F$, there is a constant $C'>0$ such that
\[r(x)\leq C'(1+|g|) \qquad\text{for }x\in g F. \]
It follows that
\[
\begin{aligned}
\int_M r(x)^{-q}\,dx
&\geq
C''\sum_{g\in\Gamma}
(1+|g|)^{-q}.
\end{aligned}
\]
The left hand side is finite because $(i+V)^{-1}\in L^q$. We conclude that $\Gamma$ has polynomial growth. Since the orbit map $\Gamma\to M$ is a quasi-isometry by the Švarc--Milnor lemma, $M$ has polynomial volume growth as well.
\end{proof}

\begin{remark}
We can rephrase Corollary \ref{lakdjnadln} as the statement that if there is a finitely summable unbounded cycle of the form $(_{f^*}L^2(M,E),D+V)$ of $C(X)\rtimes \Gamma$, with 
\begin{itemize}
\item $(_{f^*}L^2(M,E),D)$ as in Proposition \ref{adlknadlknawithf};
\item $V$ as in the assumptions of Theorem \ref{finsumthm};
\item the action of $\Gamma$ on $M$ being proper and cocompact,
\end{itemize}
then $\Gamma$ must have polynomial growth. This follows from simple the fact that for a proper, cocompact action the volume growth of the space coincides with that of the group.
\end{remark}

\section{Baum-Douglas models for equivariant $K$-homology}
\label{bdsecada}

Let $Y$ be a proper, $\Gamma$-finite $\Gamma$-CW-complex and let $\mathcal A\to Y$ be a finite-rank $\Gamma$-equivariant Azumaya bundle. We use a version of the Baum--Douglas model with coefficients in $\mathcal A$. The untwisted equivariant construction is due to Baum, Higson and Schick \cite{BaumHigsonSchick2010}; geometric models for twisted $K$-homology are developed in \cite{Barcenas2020,BaumCareyWang2013,DeeleyGoffeng2017}.

A geometric cycle for $(Y,\mathcal A)$ is a triple $(M,E,f)$ consisting of the following data:
\begin{enumerate}
\item $M$ is a proper, cocompact
$\Gamma$-spin$^{c}$ manifold without boundary;
\item $f:M\to Y$ is a proper $\Gamma$-equivariant map;
\item $E\to M$ is a finite-rank $\Gamma$-equivariant Hermitian
vector bundle together with a fiberwise non-degenerate left action $f^*\mathcal A\longrightarrow \operatorname{End}(E)$.
\end{enumerate}
The parity of the cycle is the parity of $\dim M$. We impose the equivalence relation generated by direct sum/disjoint union, $\Gamma$-spin$^{c}$ bordism and even-dimensional $\Gamma$-spin$^{c}$ vector bundle modification. The resulting $\mathbb Z/2$-graded group is denoted by $K^{\Gamma}_{*,\mathrm{geo}}(Y;\mathcal A)$.

Let $(M,E,f)$ be a cycle for $K^{\Gamma}_{*,\mathrm{geo}}(Y;\mathcal A)$. Choose a $\Gamma$-invariant complete Riemannian metric on $M$, a
$\Gamma$-invariant spin$^{c}$-connection on the spinor bundle $S_M$, and a compatible Hermitian connection on $E$. Let
\[D_E:C_c^\infty(M,E\otimes S_M) \longrightarrow C_c^\infty(M,E\otimes S_M),\]
be the corresponding twisted spin$^{c}$-Dirac operator. The left action of $f^*\mathcal A$ on $E$, together with pullback along $f$, defines a representation of $C_0(Y,\mathcal A)$ on $L^2(M,E\otimes S_M)$. Thus
\[\mu(M,E,f):=\left[f^*L^2(M,E\otimes S_M),D_E\right]\]
defines an element of $KK^\Gamma_*(C_0(Y,\mathcal A),\mathbb C)$. The usual bordism invariance of the Dirac class and the compatibility of the Thom class with vector bundle modification show that this construction descends to a homomorphism
\[\mu:K^{\Gamma}_{*,\mathrm{geo}}(Y;\mathcal A) \longrightarrow KK^\Gamma_*(C_0(Y,\mathcal A),\mathbb C).\]

\begin{thm}
\label{bhsthem}
Let $Y$ be a proper, $\Gamma$-finite $\Gamma$-CW-complex and let $\mathcal A\to Y$ be a finite-rank $\Gamma$-equivariant Azumaya bundle. Then the analytic assembly map induces an isomorphism
\[\mu:K^{\Gamma}_{*,\mathrm{geo}}(Y;\mathcal A)\longrightarrow KK^\Gamma_*(C_0(Y,\mathcal A),\mathbb C).\]
\end{thm}

\begin{proof}
For the trivial coefficient bundle this is the equivariant Baum--Douglas theorem of Baum, Higson and Schick \cite{BaumHigsonSchick2010}. The same proof applies to a finite-rank Azumaya coefficient bundle. We briefly indicate the argument for the cell-induction argument used in \cite{BaumHigsonSchick2010}; compare also the twisted equivariant construction in \cite{Barcenas2020} and the finite-rank torsion models in \cite{BaumCareyWang2013,DeeleyGoffeng2017}.

Both sides define equivariant homology theories on the category of proper $\Gamma$-CW-complexes, and the map $\mu$ is a natural transformation between them. It is therefore enough, by induction over the equivariant cells, to verify that $\mu$ is an isomorphism on an orbit $Y=\Gamma/F$, where $F\subseteq\Gamma$ is finite. Restriction to the fiber over the identity coset identifies the analytic group with $KK^F_*(\mathcal A_{eF},\mathbb C)$. Since $\mathcal A_{eF}$ is a finite-dimensional matrix algebra equipped with an $F$-action, this group is the corresponding twisted representation group in even degree and vanishes in odd degree. The geometric group has the same description by the equivariant index theorem for the finite group $F$. The index map is precisely this identification on an orbit. Induction over the
cells, using the six-term exact sequences for the analytic theory and the long exact sequences for geometric $K$-homology, proves the result.
\end{proof}

We now prove finite-summability of the geometric $K$-homology classes. Such a bound will be obtained using the isomorphism between geometric and analytic $K$-homology from the Baum-Higson-Schick theorem \ref{bhsthem}.

\begin{cor}
\label{finintesadmaodm}
If $Y$ is a proper, $\Gamma$-finite $\Gamma$-CW-complex and $\mathcal{A}\to Y$ an equivariant Azumaya bundle thereon, then any element of $KK_*^\Gamma(C_0(Y,\mathcal{A}),\C)$ can be represented by a finitely summable spectral triple. To be precise, we can take any representative to be of the form $f_*[L^2(\tilde{Y};\mathcal{E}),\slashed{D}]$ for some proper, cocompact $\Gamma$-manifold $\tilde{Y}$, equivariant mapping $f:\tilde{Y}\to Y$ and a hermitean vector bundle $\mathcal{E}\to \tilde{Y}$ of left $f^*\mathcal{A}$-modules with a Dirac operator $\slashed{D}$. If $Y$ is a proper, cocompact $\Gamma$-spin$^c$ manifold we can take $\tilde{Y}=Y$ and $f$ to be the identity. 
\end{cor}

\begin{proof}
By Theorem \ref{bhsthem}, the given class can be represented by a geometric cycle $(\widetilde Y,E,f)$. After an equivariant homotopy we may take $f$ to be locally Lipschitz with respect to a $\Gamma$-invariant piecewise Riemannian metric on $Y$. Choose $\Gamma$-invariant geometric data on $\widetilde Y$. Since the action on $\widetilde Y$ is proper and cocompact, these data have bounded geometry. Proposition \ref{adlknadlknawithf} therefore shows that $(f^*L^2(\widetilde Y,E\otimes S_{\widetilde Y}), D_E)$ is an unbounded cycle representing the given class. Proposition \ref{adlknadlkna} shows that this cycle is $p$-summable for every $p>\dim(\widetilde Y)$. Thus every element of $KK^\Gamma_*(C_0(Y,\mathcal A),\C)$ has a finitely summable unbounded representative.
\end{proof}

\section{Proof of the main results}
\label{pfsdsdsdkskn}

\subsection{Proof of Theorem \ref{mainthm1}}
\label{proofof2}

Let $\Gamma$ be a countable group with a geometric $\gamma$-element. Consider a CW-complex $X$ with a $\Gamma$-action and an element $\delta\in \gamma K^*(C(X)\rtimes \Gamma)$. Since $\gamma$ is idempotent, 
$$\delta=(1_{C(X)}\otimes \gamma)\otimes_{C(X)} \delta=(1\otimes \beta)\otimes_{C(X)\otimes \mathcal{P}}(1\otimes \alpha)\otimes_{C(X)}\delta.$$
If $X\times Z$ is a $\Gamma$-finite $\Gamma$-CW-complex, we can use Corollary \ref{finintesadmaodm} to write 
$$(1\otimes \alpha)\otimes_{C(X)}\delta=f_*[L^2(\tilde{Y};\mathcal{E}),\slashed{D}],$$ 
for some proper, cocompact $\Gamma$-manifold $\tilde{Y}$, equivariant mapping $f:\tilde{Y}\to X\times Z$ and a hermitean vector bundle $\mathcal{E}\to Y$ of left $f^*\mathcal{A}$-modules with a Dirac operator $\slashed{D}$. Therefore $\delta\in \gamma K^*(C(X)\rtimes \Gamma)$ is represented in the form
\begin{equation}
\label{lknoknadadaadaa}
\delta=(1\otimes \beta)\otimes_{C(X)\otimes A} f_*[L^2(\tilde{Y};\mathcal{E}),\slashed{D}].
\end{equation}
If $X\times Z$ is a compact spin$^c$-manifold, we can take $\tilde{Y}=X\times Z$.

Recall that by the assumption that $\gamma$ has a geometric $\beta$-factor of order $m$, $\beta$ can be represented by an $m$:th order unbounded $\Gamma$-equivariant cycle of the form $(C_0(Z,\mathcal{E}),V)$ where $\mathcal{E}\to Z$ is a hermitean vector bundle of right $\mathcal{A}$-modules and $V$ is a locally Lipschitz endomorphism valued section as in Definition \ref{geomegamadma}. It remains to prove that if $\gamma$ has a geometric $\beta$-factor of order $m\geq 1$ and $X$ is a spin$^c$-manifold, we can compute the Kasparov product \eqref{lknoknadadaadaa} as the element represented by the cycle 
\begin{equation}
\label{lknoknadadaadaacycle}
(L^2(X\times Z;\mathcal{E}\otimes E\otimes S_{X\times Z} ),1\otimes_{\nabla_\mathcal{E}}\slashed{D}+V\otimes 1_{E\otimes S_{X\times Z}}).
\end{equation}
The expression \eqref{lknoknadadaadaacycle} is an unbounded cycle for $C(X)$ of order $1$ by Proposition \ref{adlknadlknawithf}. Since $(C_0(Z,\mathcal{E}),V)$ has order $m$, the same argument shows that \eqref{lknoknadadaadaacycle} is an unbounded cycle for $C(X)\rtimes \Gamma$ of order $m$.

The technical assumption on $\beta$ will enter into the computation via the unbounded Kasparov product as recalled above in Theorem \ref{alkdnajdn}. Our computations are immediate from the approach to unbounded Kasparov products from \cite{thebeastofmesland}, or more general constructions that appeared later \cite{leschkaadlocal,leschkaadsf,leschmes}. We fix a $\Gamma$-invariant, $\mathcal{A}$-linear, hermitean connection on $E$ allowing us to lift $\slashed{D}$ to a twisted Dirac operator $1\otimes_{\nabla_E}\slashed{D}$ acting on sections of $E\otimes_{\mathcal{A}} \mathcal{E}$. We make the following observations:
\begin{enumerate}
\item By Chernoff's theorem \cite{chernoffjfa}, the twisted Dirac operator $1\otimes_{\nabla_E}\slashed{D}$ is an essentially self-adjoint operator on $L^2(X\times Z;E\otimes_{\mathcal{A}} \mathcal{E})$ because $\Gamma$ acts cocompactly and isometrically. We will tacitly identify $1\otimes_{\nabla_E}\slashed{D}$ with its self-adjoint closure in $L^2(X\times Z;E\otimes_{\mathcal{A}} \mathcal{E})$.
\item Since $V$ is locally Lipschitz, $[1\otimes_{\nabla_E}\slashed{D},V]_+$ is defined on a common core for $1\otimes_{\nabla_E}\slashed{D}$ and $V$. The anticommutator $[1\otimes_{\nabla_E}\slashed{D},V]_+$ is by the assumptions that $V^2$ is scalar valued and that $|\nabla V|\lesssim (1+|V|)^{1-\frac{1}{\alpha}}$, bounded in form sense by $C(1+|V|)^{1-1/\alpha}$ for some $C>0$. More precisely, we have that 
\begin{equation}
\label{lowered}
-C(1+|V|)^{1-1/\alpha}\leq [1\otimes_{\nabla_E}\slashed{D},V]_+\leq C(1+|V|)^{1-1/\alpha}.
\end{equation}
In particular, \cite[Theorem 1.1]{leschmes} implies that the operator 
$$1\otimes_{\nabla_E}\slashed{D}+V\otimes 1_{\mathcal{E}}$$
on $L^2(X\times Z;E\otimes_{\mathcal{A}} \mathcal{E})$ is essentially self-adjoint and its domain is given by 
\begin{align}
\label{adadklnln}
\Dom(1\otimes_{\nabla_E}\slashed{D}+&V\otimes 1_{\mathcal{E}})=\\
\nonumber
&=\Dom(1\otimes_{\nabla_E}\slashed{D})\cap (i+V)^{-1}L^2(X\times Z;E\otimes_{\mathcal{A}} \mathcal{E}\otimes S_{X\times Z}).
\end{align}
See also the results of \cite[Section 7]{leschkaadlocal}. If $\alpha=1$ and $V$ is globally Lipschitz,  it is somewhat of an overkill to invoke \cite[Theorem 1.1]{leschmes} because in this case $[1\otimes_{\nabla_E}\slashed{D},V]$ is bounded. In this case self-adjointness of the sum follows from earlier but weaker results in \cite{thebeastofmesland}.
\item Since $\Gamma$ acts cocompactly, isometrically and preserves $1\otimes_{\nabla_E}\slashed{D}$ up to zeroth order operators, we have that 
$$\Dom(1\otimes_{\nabla_E}\slashed{D})=H^1(X\times Z;E\otimes_{\mathcal{A}} \mathcal{E}\otimes S_{X\times Z}),$$
where the Sobolev space is defined from the bounded geometry arising from the $\Gamma$-action. In fact, proceeding as in \cite{boundedgeome}, we have that $(i\pm 1\otimes_{\nabla_E}\slashed{D})^{-1}$ is a pseudodifferential operator with bounded symbol of order $-1$. We can write $(i\pm (1\otimes_{\nabla_E}\slashed{D}+V))^{-1}$ as the norm-convergent limit
$$(i\pm (1\otimes_{\nabla_E}\slashed{D}+V))^{-1}=\lim_{n\to \infty} (1+\frac{i}{n}V)^{-1}(i\pm (1\otimes_{\nabla_E}\slashed{D}+V))^{-1},$$
and each limand factors over 
$$H^1(X\times Z;E\otimes_{\mathcal{A}} \mathcal{E}\otimes S_{X\times Z})\hookrightarrow L^2(X\times Z;E\otimes_{\mathcal{A}} \mathcal{E}\otimes S_{X\times Z})\xrightarrow{(1+\frac{i}{n}V)^{-1}}L^2(X\times Z;E\otimes_{\mathcal{A}} \mathcal{E}\otimes S_{X\times Z}).$$
Using that $V$ has $A$-compact resolvent, i.e. $(i\pm V)^{-1}\in C_0(Z;\End(E))$, Rellich's theorem implies that $(i\pm (1\otimes_{\nabla_E}\slashed{D}+V))^{-1}$ is compact. 
\end{enumerate}

From this discussion and Theorem \ref{alkdnajdn} we can conclude that 
\begin{equation}
\label{dadlkaldkn}
(C^\infty(X)\rtimes^{\rm alg}\Gamma, L^2(X\times Z;E\otimes_{\mathcal{A}} \mathcal{E}\otimes S_{X\times Z}),1\otimes_{\nabla_E}\slashed{D}+V\otimes 1_{\mathcal{E}}),
\end{equation}
is an $m$:th order spectral triple representing the Kasparov product \eqref{lknoknadadaadaa}. In fact, we could also use Kucerovsky's ``guess-and-check method'' directly as all conditions of \cite[Theorem 13]{dantheman}, or \cite[Theorem A.7]{GM15} for the higher order case, are satisfied by the cycle \eqref{lknoknadadaadaacycle}: item i) is a short computation, item ii) follows from Equation \eqref{adadklnln} and item iii) follows from Equation \eqref{lowered}. As such, we have the identity 
$$(1\otimes \beta)\otimes_{C(X)\otimes A} [L^2(X\times Z;\mathcal{E}\otimes S_{X\times Z}),\slashed{D}]
=[L^2(X\times Z;E\otimes_{\mathcal{A}} \mathcal{E}\otimes S_{X\times Z}),1\otimes_{\nabla_E}\slashed{D}+V\otimes 1_{\mathcal{E}}],$$
in the group 
$$KK_*^\Gamma(C(X),\C)=KK_*(C(X)\rtimes \Gamma,\C)=K^*(C(X)\rtimes \Gamma).$$

\subsection{Proof of Theorem \ref{mainthm3}}
\label{proofofthm3}

We know from Theorem \ref{mainthm1} that if $\gamma$ has a geometric $\beta$-factor of order $m$ and $X$ is a spin$^c$-manifold then any element $\delta\in \gamma K^*(C(X)\rtimes \Gamma)$ can be represented as an element of the form
$$
\delta=[L^2(X\times Z;\mathcal{E}\otimes E\otimes S_{X\times Z} ),1\otimes_{\nabla_\mathcal{E}}\slashed{D}+V\otimes 1_{E\otimes S_{X\times Z}}],$$
for some hermitean vector bundle $E\to X\times Z$ with a twisted spin$^c$-Dirac operator $\slashed{D}$ on $E\otimes S_{X\times Z}$. By Lemma \ref{lemaodknaldn}, for a finitely generated nilpotent group $\Gamma$, we can take $Z=\mathsf{N}$ to be Malcev completion of the maximal torsion-free quotient of $\Gamma$ and the $\beta$-factor to be represented by an unbounded $\Gamma$-equivariant cycle of the form $(C_0(\mathsf{N},\mathcal{E}),V)$ where $\mathcal{E}\to \mathsf{N}$ is a hermitean vector bundle of right $\mathcal{A}$-modules and $V$ is a polynomial endomorphism valued section such that $V^2$ is scalar valued and $|\nabla V|\lesssim (1+|V|)^{1-1/m}$ for some $m>0$. In particular, $(i+V)^{-1}\in L^q(\mathsf{N},\End(\mathcal{E}))$ for some $q$. By Proposition \ref{adlknadlknawithf} and Theorem \ref{finsumthm},
$$(L^2(X\times \mathsf{N};\mathcal{E}\otimes E\otimes S_{X\times \mathsf{N}} ),1\otimes_{\nabla_\mathcal{E}}\slashed{D}+V\otimes 1_{E\otimes S_{X\times \mathsf{N}}}),$$
is a finitely summable spectral triple of order $m$ on the unital $C^*$-algebra $C(X)\rtimes \Gamma$.

\section{Examples}
\label{excxca}

We turn to examples. In all the cases we consider we have geometric model for the $\gamma$-element with geometric $\beta$-factor of order $1$ with $Z$ being homeomorphic to a Euclidean space equipped with a metric depending the group.

\subsection{Dirac and Dual Dirac on $\R^n$}
We begin by describing the elements $\alpha$ and $\beta$ when $Z = \R^n$ with its flat metric and $\Gamma$ acts properly on $Z$ by isometries. It will give a geometric model for $\gamma = \beta \otimes_{C_0(\R^n,C_n)} \alpha.$ We let $C_n = \text{Cliff}(\R^n) \otimes \C,$ the complex Clifford algebra generated by the standard basis $\varepsilon_1, \ldots, \varepsilon_n$ of $\R^n$ under the relation $\varepsilon_i\varepsilon_j + \varepsilon_j \varepsilon_i = 2 \delta_{ij}$ and $\varepsilon_i^* = \varepsilon_i.$ Then $C_0(\R^n,C_n)$ is the $C^*$-algebra of continuous sections of the Clifford bundle that vanish at infinity that carries the $\Gamma$-action induced by the $\Gamma$ action on $\R^n.$ 

\begin{prop}
\label{RnDD}
    The element $\beta \in KK^\Gamma(\C,C_0(\R^n, C_n))$ is represented by the cycle $(C_0(\R^n, C_n),V)$ where $C_0(\R^n, C_n)$ is a Hilbert module over itself with the $\Gamma$-action induced by the action on $\R^n$ and $V$ given by multiplication by the function $V: \R^n \to C_n$ where $V(x) =x.$
\end{prop}

This result follows from \cite{kaspinv}. Observe that $V$ is Lipschitz with constant 1 and, further, that $V$ satisfies the conditions of being a geometric $\beta$ factor as in Definition \ref{geomegamadma} as, selecting $x_0$ as the origin: $|V(x)| = |x| = d(x,x_0)$. The next result can also be found in \cite{kaspinv}.

\begin{prop}
\label{RnD}
The class $\alpha \in KK^\Gamma(C_0(\R^n,C_n), \C)$ is represented by the cycle 
$$(L^2(\R^n,\wedge^*\C^n),D),$$ 
where the Hilbert space $L^2(\R^n,\wedge^*\C^n)$ carries the $C_0(\R^n, C_n)$ representation via pointwise Clifford multiplication. The action of $\Gamma$ on $L^2(\R^n,\wedge^*\C^n)$ is induced by the action of $\Gamma$ on $\R^n.$ Then we define the Dirac operator $D$ on $L^2(\R^n,\wedge^*\C^n)$ as the closure of the differential operator
$$D = \displaystyle - i \sum_j \varepsilon_j \frac{\partial}{\partial x_j},$$ 
with domain $C^\infty_c(\R^n,\wedge^*\C^n)$.
\end{prop} 

A short computation with harmonic oscillators show that $\alpha$ and $\beta$ as constructed above (when $\Gamma =\{*\})$, satisfies 
$$\gamma = \beta \otimes_{C_0(\R^n, C_n)} \alpha = 1 \in KK(\C,\C).$$
More generally, the following holds.

\begin{prop}{\cite{Val}}
 When $\Gamma$ is a countable group acting properly on $\R^n$ by isometries, $\gamma = \beta \otimes_{C_0(\R^n, C_n)} \alpha = 1 \in KK^\Gamma(\C,\C).$
\end{prop}

In the context of Theorem \ref{mainthm1}, this means that, if $\Gamma$ is a countable group acting properly on $\R^n$ by isometries and $X$ is a $\text{spin}^c$-manifold, then every class in $K^*(C(X) \rtimes \Gamma)$ can be represented by a cycle of the form $(L^2(X \times \R^n, C_n \otimes_{C_n} \mathcal{E}), 1 \otimes_{\nabla_{C_n}} D + V \otimes 1_{C_n})$, where $\mathcal{E}$ is a Clifford bundle on $X$ with Dirac operator $D$. Further, $1 \otimes \nabla_{C_n} D$ reduces to $D$ under the isomorphism $C_n \otimes_{C_n} \mathcal{E} \cong \mathcal{E}$. Since $(1+|x|)^{-1}\in L^p(\R^n)$ for any $p> n$, Theorem \ref{finsumthm} ensures that the cycles $(L^2(X \times \R^n, \mathcal{E}), D + V)$ are finitely summable.

\subsection{Fundamental Groups of Flat Manifolds}
Let $M$ be a compact Riemannian flat manifold. Then $M \cong \R^n / \pi_1(M),$ where $\pi_1(M)$ is a \textit{Bieberbach group} acting on $\R^n$ by rigid motions \cite{Char}. i.e. an element $g \in \pi_1(M)$ can be written as $(m,s)$ where $m \in O(n), s \in \R^n,$ so that $m$ corresponds to rotation on $\R^n$ and $s$ to translation. Elements are composed via $(m_1,s_1)(m_2,s_2) = (m_1m_2, m_1s_2 + s_1)$.

\begin{prop}
\label{propbieber}
Let $M$ be a Riemannian flat manifold and $\Gamma \cong \pi_1(M).$ If $\Gamma$ acts on a compact $\text{spin}^c$-manifold $X,$ each element of $K^*(C(X) \rtimes \Gamma)$ can be represented by a finitely summable unbounded cycle of the form $L^2(X \times \R^n,E), D +V')$ where $E$ is a Clifford bundle with Dirac operator $D$ and, $V'(f)(x,z) = zf(x,z),$ for $f \in L^2(X \times \R^n, E), x\in X, z \in \R^n.$ 
\end{prop}

\begin{proof}
Bieberbach groups are torsion-free and discrete. Further, $M$ is universally covered by $\R^n.$ Thus, we can take $Z,$ as in Definition \ref{geomegamadma}, to be $\R^n.$ If $\Gamma \cong \pi_1(M)$ for a flat manifold, $M,$ $\pi_1(M)$ acts on $\R^n$ and we obtain $\alpha$ and $\beta$ as $KK^\Gamma(C_0(\R^n, C_n), \C)$ and $KK^\Gamma(\C, C_0(\R^n, C_n))$ cycles as in Propositions \ref{RnDD} and \ref{RnD} because $\Gamma$ acts properly on $\R^n$ by isometries. Therefore, $\gamma =1 \in KK^{\Gamma}(\C,\C)$ has a geometric $\beta$ factor. Then, if $\Gamma$ acts on a compact $\text{spin}^c$-manifold $X,$ the diagonal action of $\Gamma$ on $X \times \R^n$ is proper. The $\Gamma$-action on $X \times \R^n$ is cocompact because $M=\R^n/\Gamma$ is compact. 

We can from Corollary \ref{finintesadmaodm} conclude that each class in $KK^\Gamma(C_0(X \times \R^n, C_n), \C)$ can be represented by a cycle of the form $(L^2(X \times Z, E), D)$ for a hermitian vector bundle $E$ with Dirac operator $D.$ 
 Hence, for $\delta\in \gamma K^*(C(X) \rtimes \Gamma) \cong K^*(C(X) \rtimes \Gamma)$, using $\gamma =1_\mathbb{C}$, we have that: 
 \begin{align*}
    \delta &= (1_{C(X)} \otimes \beta) \otimes_{C(X) \otimes C_0(\R^n,C_1)} [L^2(X\times \R^n; \mathcal{E},D)] \\ &= [C_0(X\times \R^n, C_n), 1_{X} V'] \otimes_{C_0(X\times \R^n, C_n)} [L^2(X \times \R^n), D] \\ &= [L^2(X\times \R^n; \mathcal{E}), D + V']
\end{align*} 
where $V'$ acts on $L^2(X \times \R^n, \mathcal{E})$ by (clifford) multiplication by the function $f(x,z) = z$ (where $x \in X, z \in \R^n),$ noting the identification 
\[C_0(X \times \R^n, C_n) \otimes_{C_0(X \times \R^n)} L^2(X \times \R^n, C_n \otimes_{C_n} \mathcal{E}) \cong L^2(X\times \R^n, C_n \otimes_{C_n} E) \cong L^2(X \times \R^n, \mathcal{E}).\] 
Note that the $\Gamma$-action on $L^2(X \times \R^n, C_n \otimes_{C_n} \mathcal{E})$ is induced from the diagonal action on $X \times \R^n$ and a function $f \in C(X)$ acts by Clifford multiplication via $f\cdot \xi (x,z) = f(x)\xi(x,z)$. The Kasparov product with the class $1 \otimes \beta \in KK^\Gamma(C(X),C_0(X \times \R^n, C_n))$ is by Theorem \ref{alkdnajdn} the class $[L^2(X \times Z,E), D+V' ]$. Further, $(i+V)^{-1} \in L^p(\R^n,C_n)$ for all $p > n,$ so these cycles are finitely summable as well by Theorem \ref{finsumthm}.
\end{proof}

\begin{example}[Irrational rotation algebras]
Consider the group $\Gamma = \Z$ acting on $Z = \R$ via translation. Then $\Z$ has Dirac element as in Proposition $\ref{RnD}$ and geometric $\beta$ factor in $KK^\Z(\C,C_0(\R,C_1))$ as in Proposition \ref{RnDD}. As in \cite{Ech}, because $\Z$ is a closed subgroup of $\R,$ we have that $\gamma =\beta \otimes_{C_0(\R,C_1)} \alpha =1$. For a $\Z$-space $X$, the quotient space $X \times \R/ \: \Z$, formed from the diagonal $\Z$ action, is the \textit{mapping torus} of the $\Z$.action on $X$.
 
Consider $X = S^1$ with an action of $\Gamma = \Z$ by rotation with the angle $2\pi \alpha$. The $C^*$-algebra $C(S^1) \rtimes_\alpha \Z$ is usually called the irrational rotation algebra, or the noncommutative torus. In this case $X \times Z \cong S^1 \times \R.$ Thus, for any $\delta \in K^*(C(S^1) \rtimes_\alpha \Z),$ \[(1_{C(S^1)} \otimes \beta) \otimes_{C(S^1 \times \R, C_1)} \delta\in KK^\Z(C_0(S^1 \times \R, C_1),\C)\] where $\Z$ acts diagonally on $S^1 \times \R.$ Thus, $S^1 \times \R/ \: \Z$ is the mapping torus of the irrational rotation action. Then, using Corollary \ref{finintesadmaodm}, any class in $KK^\Z(C_0(S^1 \times \R, C_1), \C)$ can be represented by a cycle of the form $(L^2(S^1 \times \R, \mathcal{E}), \mathcal{D}),$ where $\mathcal{D}$ is $p > 2$ summable. By Proposition \ref{propbieber}, for $\delta\in K^*(C(S^1) \rtimes \Z),$ we have that: 
   \[\delta =[L^2(S^1 \times \R; \mathcal{E}), D + V']\]
 where $V'$ acts on $L^2(S^1 \times \R, \mathcal{E})$ by Clifford multiplication by the function $f(x,z) = z$, for $x \in S^1$, $z \in \R$. Since $(i +V)^{-1} \in L^p(\R,C_1)$ for all $p >1$,  these cycles are finitely summable. 
\end{example}

The same construction as in the examples above applies to the Klein bottle $K$ with $\Gamma=\pi_1(K)$ acting on the double cover $\mathbb{T}^2$ of $K$, and for the ten compact flat $3$-manifolds \cite{Wolf} where we consider fundamental groups acting on $\mathbb{R}^3$ by deck transformations. .

\subsection{Proper special $\Gamma$-manifolds}

We have the following more geometric formulation of Theorem \ref{mainthm3} when $\Gamma$ acts on a geometry of nonpositive sectional curvature with polynomial growth. The polynomial growth assumption is crucial for our finite-summability proof to go through.

\begin{prop}
Let $M$ be a $\Gamma$-compact, proper $\Gamma$-manifold with nonpositive sectional curvature and polynomial volume growth. Fix $x \in M$ and denote the geodesic distance function on $M$ by $\rho$. If $\Gamma$ acts on a compact manifold $X$, each class in $K^*(C(X) \rtimes \Gamma)$ can be represented by a finitely summable unbounded cycle of the form 
\[(L^2(X \times M, \mathcal{E}), D + \delta')\] 
where $\mathcal{E}\to X \times M$ is a Clifford bundle with Dirac operator $D$ and $\delta'$ is Clifford multiplication by the one-form $\delta'(x,z) = \rho(x_0,z)d_x \rho(x_0,z).$
\end{prop}

\begin{proof}
As noted in Remark \ref{rempos}, if $M$ is a proper $\Gamma$-manifold with nonpositive sectional curvature it is a special $\Gamma$-manifold as in Definition \ref{special}. In fact Proposition \ref{betarep} guarantees that the $\gamma$-element has the geometric $\beta$-factor $\beta \in KK^\Gamma(\C, C_0(M, \C\ell(M))$  given by $(C_0(M, \C\ell(M)),\delta)$ by Proposition \ref{betarep}. So, when $\Gamma$ acts on a compact $X$, we have that $1 \otimes_{C(X)} \beta = (C_0(X \times M, \C \ell(M)), \delta'),$ where $\delta'$ is multiplication by the one-form $\delta'(x,z) = \rho(x_0,x)d_x \rho(x_0,x).$ Further, the diagonal action of $\Gamma$ on $X \times M$ is proper. Then, if $(L^2(X \times M,\mathcal{E}), D)$ is a cycle giving an element of $KK^\Gamma(C_0(X \times M, \C \ell(M)),\C)$, 
\[(1\otimes_{C(X)} \beta) \otimes_{C(X\times M,\mathcal{E})} [L^2(X \times M,\mathcal{E}), D] = [L^2(X \times M, \mathcal{E}), D + \delta'] .\] 
 By the polynomial volume growth of $M$, $(i+\delta')^{-1} \in L^p(X \times M, \End(\mathcal{E}))$ for some $p$ so that, by Theorem \ref{finsumthm}, this cycle is finitely summable. 
 
By Theorem \ref{mainthm1}, the classes defined from cycles of the form $(L^2(X \times M, \mathcal{E}), D + \delta')$ exhaust  $\gamma K^*(C(X) \rtimes \Gamma)$. Since $\Gamma$ acts properly on $M$, and the latter has polynomial growth, then so does $\Gamma$. In particular, $\Gamma$ satisfies the strong Baum-Connes conjecture and $K^*(C(X) \rtimes \Gamma)=\gamma K^*(C(X) \rtimes \Gamma)$. 
\end{proof}

\subsection{Discrete subgroups in semisimple Lie groups}

In the other extreme compared to the above examples of groups acting isometrically on flat $\R^n$ we have discrete subgroups in semisimple Lie groups. Assume that $G$ is a connected, semisimple Lie group, and write $K\subseteq G$ for a choice of maximal compact subgroup. Whenever $\Gamma\subseteq G$ is a discrete subgroup, $\Gamma$ acts properly on the symmetric space $G/K$. Since $G/K$ has nonpositive sectional curvature, there is a diffeomorphism $G/K\cong \R^n$. This diffeomorphism is explicit in terms of the Cartan decomposition and it is not isometric for the flat metric on $\R^n$, rather $G/K$ has exponential volume growth. In particular,  the above finite-summability arguments all collapse for $Z=G/K$. By Puschnigg's work \cite{Pusch}, this is to be expected as lattices in higher rank semisimple Lie groups have no finitely summable Fredholm modules.

We specialize to the case when $X=G/P$ for a parabolic subgroup $P\subseteq G$ and $\Gamma\subseteq G$ is torsion free. We let $\Gamma$ act on the compact manifold $X:=G/P$. Then 
$$Z=\underline{\mathcal{E}}\Gamma=G/K \quad\mbox{and}\quad X\times Z\cong G/(P\cap K).$$
The last $G$-equivariant diffeomorphism is defined by 
$$G/(P\cap K)\to G/P\times G/K,\quad g(P\cap K)\mapsto (gP,gK),$$ 
which is surjective for rank reasons and injective by equivariance. If $\Gamma$ is torsion-free and cocompact, $E\Gamma=G/K$ and 
$$X\times E\Gamma=G/(P\cap K)\to (X\times E\Gamma)/\Gamma=\Gamma \backslash G/(P\cap K),$$
with the last space forming a fiber bundle over the locally symmetric space $B\Gamma=\Gamma\backslash G/K$ with fiber $K/(P\cap K)=G/P=X$.

\begin{thm}
Any class in $\gamma K^*(C(G/P)\rtimes \Gamma)$ can be represented by a cycle of the form $(L^2(G/(P\cap K),S_{G/(P\cap K)}\otimes E, D_E+V)$, where 
$V(g(P\cap K))=\delta(gK)$, 
for $\delta\in C^\infty(G/K,\C\ell(G/K))$ defined by
$$\delta(x):=\rho(x)\rd\rho(x),$$
where $\rho(x)=\rho(eK,x)$ for the geodesic distance $\rho$ on $G/K$. These cycles are never finitely summable.
\end{thm}

\begin{proof}
The result follows from Theorem \ref{mainthm1} and Proposition \ref{betarep}. We note that in Theorem \ref{mainthm1} we need the technical assumption that $\Gamma$ acts cocompactly to ensure bounded geometry and compact resolvent of the Kasparov product, however this is clear for any discrete subgroup $\Gamma$ in the case at hand. The symmetric space $G/K$ always has exponential growth, so the cycles in question are never finitely summable by Corollary \ref{lakdjnadln}.
\end{proof}

\begin{remark}
In the special case that $P$ is the minimal parabolic, we can remove the assumption on the class to belong to the $\gamma$-part of $K$-homology. Indeed, we have that
$$\gamma K^*(C(G/P)\rtimes \Gamma)= K^*(C(G/P)\rtimes \Gamma),$$
when $P$ is the minimal parabolic. This follows from the following argument. The minimal parabolic can be Langlands decomposed as $P=MAN$ where $M$ is compact and $AN$ is solvable, so $P$ is amenable. In particular, $P$ satisfies the strong Baum-Connes conjecture, i.e. $\gamma_P=1$. Therefore $\C$ is compactly induced in $KK^P$ in the Meyer-Nest sense \cite{meyernest}, and so $C(G/P)=\mathrm{Ind}_P^G\C$ is compactly induced in $KK^G$ in the Meyer-Nest sense. We therefore have that $1_{C(G/P)}\otimes \gamma_G=1_{C(G/P)}$ so $\gamma_G$ acts as the identity on $C(G/P)$ in $KK^G$. We conclude that $\gamma_\Gamma$ acts as the identity on $C(G/P)$ in $KK^\Gamma$.
\end{remark}

\end{document}